\documentclass[a4paper,12pt,twoside]{article}
\usepackage[francais]{babel}
\usepackage{amssymb, amsmath, eucal}
\usepackage[all]{xy}
\usepackage{mathrsfs}
\usepackage[dvips]{graphics}
\begin{document}
\textwidth15.5cm
\textheight22.5cm
\voffset=-13mm
\newtheorem{The}{Theorem}[section]
\newtheorem{Lem}[The]{Lemma}
\newtheorem{Prop}[The]{Proposition}
\newtheorem{Cor}[The]{Corollary}
\newtheorem{Rem}[The]{Remark}
\newtheorem{Obs}[The]{Observation}
\newtheorem{SConj}[The]{Standard Conjecture}
\newtheorem{Titre}[The]{\!\!\!\! }
\newtheorem{Conj}[The]{Conjecture}
\newtheorem{Question}[The]{Question}
\newtheorem{Prob}[The]{Problem}
\newtheorem{Def}[The]{Definition}
\newtheorem{Not}[The]{Notation}
\newtheorem{Claim}[The]{Claim}
\newtheorem{Conc}[The]{Conclusion}
\newtheorem{Ex}[The]{Example}
\newtheorem{Fact}[The]{Fact}
\newcommand{\C}{\mathbb{C}}
\newcommand{\R}{\mathbb{R}}
\newcommand{\N}{\mathbb{N}}
\newcommand{\Z}{\mathbb{Z}}
\newcommand{\Q}{\mathbb{Q}}
\newcommand{\Proj}{\mathbb{P}}

\begin{center}

{\Large\bf Compact Complex Manifolds with Small Gauduchon Cone}

\end{center}

\begin{center}

{\large Dan Popovici and Luis Ugarte}

\end{center}

\vspace{1ex}

\noindent {\small {\bf Abstract.} This paper is intended as the first step of a programme aiming to prove in the long run the long-conjectured closedness under holomorphic deformations of compact complex manifolds that are bimeromorphically equivalent to compact K\"ahler manifolds, known as Fujiki {\it class} ${\cal C}$ manifolds. Our main idea is to explore the link between the {\it class} ${\cal C}$ property and the closed positive currents of bidegree $(1,\,1)$ that the manifold supports, a fact leading to the study of semi-continuity properties under deformations of the complex structure of the dual cones of cohomology classes of such currents and of Gauduchon metrics. Our main finding is a new class of compact complex, possibly non-K\"ahler, manifolds defined by the condition that every Gauduchon metric be strongly Gauduchon (sG), or equivalently that the Gauduchon cone be small in a certain sense. We term them sGG manifolds and find numerical characterisations of them in terms of certain relations between various cohomology theories (De Rham, Dolbeault, Bott-Chern, Aeppli). We also produce several concrete examples of nilmanifolds demonstrating the differences between the sGG class and well-established classes of complex manifolds. We conclude that sGG manifolds enjoy good stability properties under deformations and modifications.}

\vspace{1ex}

\noindent Mathematics subject classification (2010)\!: 32G05, 53C55, 14C30, 14F40.

\section {Introduction}

Fujiki {\it class} ${\cal C}$ manifolds provide a key link between K\"ahler and non-K\"ahler geometries that justifies the central role they play in the classification theory of compact complex manifolds. They are defined (cf. [Fuj78], [Var86]) by the condition that for any such compact complex manifold $X$ there exist a holomorphic bimeromorphic map ($=$ a modification) $\mu:\widetilde{X}\longrightarrow X$ from a compact K\"ahler manifold $\widetilde{X}$. Deep insights into the classification theory are often provided by the theory of deformations of complex structures.

  It has long been conjectured that the deformation limit of any holomorphic family of {\it class} ${\cal C}$ manifolds ought to be a {\it class} ${\cal C}$ manifold.

\begin{Conj}\label{Conj:classC_limit} Let $\pi:{\cal X}\to\Delta$ be a proper holomorphic submersion between a complex manifold ${\cal X}$ and an open disc $\Delta\subset\C$ containing the origin. Suppose that the fibre $X_t:=\pi^{-1}(t)$ is a {\it class} ${\cal C}$ manifold for every $t\in\Delta\setminus\{0\}$. 

Then the remaining (limit) fibre $X_0:=\pi^{-1}(0)$ is again a {\it class} ${\cal C}$ manifold.

\end{Conj}

 A two-step strategy for tackling this conjecture was briefly outlined in [Pop13b]\!\!: \\

{\it Step} $1$\!\!: {\it prove that a compact complex manifold $X$ belongs to the class ${\cal C}$ if and only if there are ``many'' closed positive $(1,\,1)$-currents on $X$.}

\vspace{2ex}

{\it Step} $2$\!\!: {\it prove that there can only be ``more'' closed positive $(1,\,1)$-currents on $X_0$ than on the generic fibre $X_t$.}

\vspace{2ex}

 We now make the meaning of these goals explicit. Our approach relies on the interaction between two points of view\!\!: cohomology and positivity.

\vspace{2ex}

\noindent {\it Reminder of a few definitions}

\vspace{2ex}

Throughout the paper, $X$ will stand for a possibly non-K\"ahler compact complex manifold of complex dimension $n$.

The following two cohomologies are especially relevant in the possibly non-K\"ahler context. For $p,q=0,\dots , n$, the {\bf Bott-Chern}, resp. {\bf Aeppli cohomology groups} of bidegree $(p,\,q)$ are defined as

$$H^{p,\,q}_{BC}(X,\,\C)= \frac{\ker\partial\cap\ker\bar\partial}{\mbox{Im}\,\partial\bar\partial}, \hspace{3ex} H^{p,\,q}_A(X,\,\C)= \frac{\ker\partial\bar\partial}{\mbox{Im}\,\partial + \mbox{Im}\,\bar\partial},$$  

\noindent where all the vector spaces involved are subspaces of the space $C^{\infty}_{p,\,q}(X,\,\C)$ of smooth $(p,\,q)$-forms on $X$. Both of these types of cohomology groups, much like all the other familiar cohomologies, can be computed using either smooth forms or currents.

 On the positivity side, since divisors may not exist on non-K\"ahler manifolds, one turns to their transcendental analogues that are provided by closed currents of bidegree $(1,\,1)$. One of our main tools will be the by now classical {\bf pseudo-effective cone} of a compact complex manifold $X$ (cf. [Dem92])\!\!:

$${\cal E}_X:=\{[T]_{BC}\in H^{1,\,1}_{BC}(X,\,\R)\,\mid\,T\hspace{1ex}\mbox{is a closed positive}\hspace{1ex}(1,\,1)-\mbox{current on}\hspace{1ex}X\}.$$

\noindent This is a closed convex cone in $H^{1,\,1}_{BC}(X,\,\R)$.

 On the other hand, since we do not assume the existence of K\"ahler metrics on the manifolds that we shall study below, Gauduchon metrics will play a key role. Recall that a $C^{\infty}$ positive definite $(1,\,1)$-form $\omega>0$ on an $n$-dimensional $X$ is said to be a {\bf Gauduchon metric} ([Gau77a]) if $\partial\bar\partial\omega^{n-1}=0$, while $\omega$ is a {\bf strongly Gauduchon (sG) metric} ([Pop13a]) if $\partial\omega^{n-1}\in\mbox{Im}\,\bar\partial$. The main virtue of Gauduchon metrics is that they always exist ([Gau77a]), providing a useful substitute for K\"ahler metrics when the latter do not exist.

 Another key tool in this paper will be the {\bf Gauduchon cone} of $X$ (cf. [Pop13b]):

$${\cal G}_X:=\{[\omega^{n-1}]_A\in H^{n-1,\,n-1}_A(X,\,\R)\,\mid\,\omega\hspace{1ex}\mbox{is a Gauduchon metric}\hspace{1ex} \mbox{on}\hspace{1ex}X\}.$$

\noindent This is an open convex cone in $H^{n-1,\,n-1}_A(X,\,\R)$. Also recall that the {\bf sG cone} was defined in [Pop13b, section $\S.5$] as

$${\cal SG}_X = {\cal G}_X\cap\ker T\subset{\cal G}_X\subset H^{n-1,\,n-1}_A(X,\,\R),$$

\noindent i.e. the intersection of the Gauduchon cone with the kernel of the following canonical linear map

\begin{equation}\label{eqn:T-def}T\,:\,H^{n-1,\,n-1}_A(X,\,\C)\longrightarrow H^{n,\,n-1}_{\bar\partial}(X,\,\C), \hspace{3ex} T([\Omega]_A):=[\partial\Omega]_{\bar\partial},\end{equation}

\noindent for any $[\Omega]_A\in H^{n-1,\,n-1}_A(X,\,\C)$. The map $T$ is well defined (i.e. independent of the choice of representative of the Aeppli class $[\Omega]_A$) and shows that the sG property is cohomological\!\!: either all the Gauduchon metrics $\omega$ for which $\omega^{n-1}$ belongs to a given Aeppli-Gauduchon class $[\omega^{n-1}]_A\in{\cal G}_X$ are strongly Gauduchon (in which case we say that $[\omega^{n-1}]_A$ is an sG class), or none of them is. In other words, the sG cone ${\cal SG}_X$ is the set of all sG classes on $X$. It is empty if $X$ does not support any sG metric.

 Moreover, under the duality (see e.g. [Sch07] or [Pop13b])

\begin{equation}\label{eqn:BC-A-duality} H^{1,\,1}_{BC}(X,\,\C)\times H^{n-1,\,n-1}_A(X,\,\C)\longrightarrow\C, \hspace{3ex} ([\alpha]_{BC},\,[\beta]_A)\mapsto\int\limits_X\alpha\wedge\beta,\end{equation}

\noindent the pseudo-effective cone ${\cal E}_X\subset H^{1,\,1}_{BC}(X,\,\R)$ is dual to the closure of the Gauduchon cone $\overline{{\cal G}_X}\subset H^{n-1,\,n-1}_A(X,\,\R)$ of $X$ thanks to Lamari's duality lemma (cf. [Lam99, Lemme 3.3]).

\vspace{2ex}

\noindent {\it Back to the two-step approach to Conjecture \ref{Conj:classC_limit}}

\vspace{2ex}

 {\it Step} $1$ would be the transcendental analogue of the following well-known fact\!\!: a compact complex manifold $X$ is Moishezon ($=$ bimeromorphically equivalent to a projective manifold) if and only if there are ``many'' divisors on $X$ (in the sense that the algebraic dimension of $X$ is maximal, i.e. equal to the dimension of $X$ as a complex manifold).

 The meaning of ``many'' in connection with closed positive $(1,\,1)$-currents has yet to be probed, but we suspect that it will mean that the pseudo-effective cone ${\cal E}_X$ of $X$ is ``maximal'' at least in the following sense\!\!:

$$(i)\,\,\mathring{{\cal E}}_X\neq\emptyset \hspace{2ex} \mbox{and} \hspace{2ex} (ii)\,\,{\cal SG}_X = {\cal G}_X,  \hspace{10ex}(\star)$$

\noindent where $\mathring{}$ stands for ``interior'', while ${\cal SG}_X$ and ${\cal G}_X$ are respectively the sG cone and the Gauduchon cone of $X$. Property $(i)$ uses the non-emptiness of the interior as a way of requiring ${\cal E}_X$ to be fairly large, while property $(ii)$ requires ${\cal G}_X$ to be fairly small, hence by duality ${\cal E}_X$ to be again fairly large by a different criterion.

 Each of the two properties in $(\star)$ is necessary for $X$ to be of class ${\cal C}$\footnote{The interior of the pseudo-effective cone ${\cal E}_X$ contains the big cone of cohomology classes of K\"ahler currents (cf. e.g. [BDPP, $\S.3$]), while the latter is non-empty (and equals $\mathring{{\cal E}}_X$) on {\it class} ${\cal C}$ manifolds by [DP04, Theorem 3.4]. Meanwhile, any {\it class} ${\cal C}$ manifold satisfies the $\partial\bar\partial$-lemma by [AB95], which, in turn, implies $(ii)$ in $(\star)$ by Lemma \ref{Lem:sGG-equiv-descriptions1}.}, but none of them is sufficient on its own (see section \ref{section:interior-psef} for examples of manifolds not in the class ${\cal C}$ whose pseudo-effective cone has non-empty interior). However, together they may become sufficient, or should condition $(\star)$ turn out to be insufficient for $X$ to be in the {\it class} ${\cal C}$, it will have to be reinforced.

\vspace{2ex}

 {\it Step} $2$ means that the pseudo-effective cone ${\cal E}_{X_t}$ can only increase in the limit as $t\to 0$ (i.e. it behaves upper-semicontinuously under deformations of the complex structure of $X_t$), while its dual, the (closure of the) Gauduchon cone ${\cal G}_{X_t}$, can only decrease in the limit (i.e. it behaves lower-semicontinuously).

\vspace{2ex}

 In this paper, we begin the implementation of this two-step strategy by studying the manifolds defined by property $(ii)$ in $(\star)$ and by giving a complete affirmative answer to the problem raised at {\it Step} $2$ of this line of argument for deformations of such manifolds.

\vspace{2ex}

\noindent {\it Definitions introduced and results obtained in this paper}

\vspace{2ex}

 $(1)$\, In line with the goals in the first step of the approach to Conjecture \ref{Conj:classC_limit} outlined above, we set out to investigate the following class of compact complex manifolds introduced in [Pop13b] for which we now propose the following terminology.

\begin{Def}\label{Def:sGG} Let $X$ be a compact complex manifold, $\mbox{dim}_{\C}X=n$.  We say that $X$ is an {\bf sGG manifold} if the sG cone of $X$ coincides with the Gauduchon cone of $X$, i.e. if ${\cal SG}_X = {\cal G}_X$.

\end{Def}

  Since the kernel of the linear map $T$ defined in (\ref{eqn:T-def}) is a vector subspace of $H^{n-1,\,n-1}_A(X,\,\C)$, its intersection with the open convex Gauduchon cone leaves the latter unchanged if and only if $\ker T = H^{n-1,\,n-1}_A(X,\,\C)$, i.e. if and only if $T$ vanishes identically. We obtain equivalent descriptions of the sGG property summed up as follows (cf. [Pop13b, section $\S.5$] for $(i)-(iii)$).

\begin{Lem}\label{Lem:sGG-equiv-descriptions1} The following statements are equivalent\!\!:

\vspace{1ex}

\noindent $(i)$\, $X$ is an sGG manifold\!;

\vspace{1ex}

\noindent $(ii)$\, the map $T$ vanishes identically\!;

\vspace{1ex}

\noindent $(iii)$\, the following special case of the $\partial\bar\partial$-lemma holds\!\!: for every $d$-closed $(n,\,n-1)$-form $\Gamma$ on $X$, if $\Gamma$ is $\partial$-exact, then $\Gamma$ is also $\bar\partial$-exact\!;

\vspace{1ex}

\noindent $(iv)$\, every Gauduchon metric $\omega$ on $X$ is strongly Gauduchon.

\end{Lem}

 Recall that a compact complex manifold $X$ is said to be a $\partial\bar\partial$-manifold if for every bidegree $(p,\,q)$ and every smooth $d$-closed $(p,\,q)$-form $u$ on $X$, the $\partial$-exactness, the $\bar\partial$-exactness, the $d$-exactness and the $\partial\bar\partial$-exactness of $u$ are pairwise equivalent. The $\partial\bar\partial$-property is equivalent to all the canonical morphisms $H^{p,\,q}_{BC}(X,\,\C)\to H^{p,\,q}_A(X,\,\C)$ being isomorphisms and implies the Hodge decomposition and the Hodge symmetry, as well as the degeneration at $E_1$ of the Fr\"olicher spectral sequence of $X$. On the other hand, every {\it class} ${\cal C}$ manifold is a $\partial\bar\partial$-manifold,  but there exist $\partial\bar\partial$-manifolds that are not of  {\it class} ${\cal C}$. (See e.g. [Pop14] for a review of these matters.) Lemma \ref{Lem:sGG-equiv-descriptions1} shows, in particular, that every $\partial\bar\partial$-manifold is sGG. We shall see below (cf. e.g. Corollary \ref{Cor:Iwasawa-sGG}) that the converse fails, so the sGG class strictly contains the $\partial\bar\partial$ class.

 We start by giving two more characterisations of sGG manifolds. The first one is a numerical characterisation in terms of the Bott-Chern number $h^{0,\,1}_{BC}:=\mbox{dim}_{\C}H^{0,\,1}_{BC}(X,\,\C)$ and the Hodge number $h^{0,\,1}_{\bar\partial}:=\mbox{dim}_{\C}H^{0,\,1}_{\bar\partial}(X,\,\C)$.

\begin{The}\label{The:BC-dbar} On any compact complex manifold $X$ we have $h^{0,\,1}_{BC}\leq h^{0,\,1}_{\bar\partial}$.

 Moreover, $X$ is an sGG manifold if and only if $h^{0,\,1}_{BC} = h^{0,\,1}_{\bar\partial}$.

\end{The}

An immediate consequence is the following

\begin{Cor}\label{Cor:Iwasawa-sGG} The Iwasawa manifold and all its small deformations in its Kuranishi family are sGG manifolds (but, of course, not $\partial\bar\partial$-manifolds).

\end{Cor}

Recall that the Iwasawa manifold is the nilmanifold of complex dimension $3$ obtained as the quotient of the Heisenberg group of $3\times 3$ upper triangular matrices with entries in $\C$ by the subgroup of those matrices with entries in $\Z[i]$. Historically, it was the first compact complex manifold to be discovered whose Fr\"olicher spectral sequence does not degenerate at $E_1$. In particular, it is not a $\partial\bar\partial$-manifold. Thanks to Corollary \ref{Cor:Iwasawa-sGG}, the Iwasawa manifold is our main example of sGG manifold that is not $\partial\bar\partial$. Its Kuranishi family was explicitly computed by Nakamura in [Nak75].

\vspace{2ex}

 Our second numerical characterisation of sGG manifolds is in terms of the first Betti number $b_1:=\mbox{dim}_{\C}H^1_{DR}(X,\,\C)$ and the Hodge number $h^{0,\,1}_{\bar\partial}$.

\begin{The}\label{The:Betti-dbar} On any compact complex manifold $X$ we have $b_1\leq 2h^{0,\,1}_{\bar\partial}$.

 Moreover, $X$ is an sGG manifold if and only if $b_1 = 2h^{0,\,1}_{\bar\partial}$.

\end{The}

 This makes sGG manifolds reminiscent of compact K\"ahler surfaces\!\!: recall that a compact complex surface is K\"ahler if and only if $b_1$ is even (cf. Kodaira's classification of surfaces, [Miy74] and [Siu83], or [Buc99] and [Lam99] for a direct proof). In dimension $2$, the K\"ahler and sGG conditions are clearly equivalent, but in dimension $\geq 3$ the sGG property is tremendously weaker than the K\"ahler one.

 We infer from Theorem \ref{The:Betti-dbar} that the sGG property of compact complex manifolds is {\bf open} under holomorphic deformations $(X_t)_{t\in\Delta}$. We denote by $h^{p,\,q}_{\bar\partial}(t)$, resp. $h^{p,\,q}_{BC}(t)$, the Hodge number, resp. Bott-Chern number, of bidegree $(p,\,q)$ of the fibre $X_t$ for any $t$.

\begin{Cor}\label{Cor:openness} Let $(X_t)_{t\in\Delta}$ be any holomorphic family of compact complex manifolds. Fix an arbitrary $t_0\in\Delta$. If $X_{t_0}$ is an sGG manifold, then\!\!:

\vspace{1ex}

\noindent $(i)$\, $X_t$ is an sGG manifold for all $t\in\Delta$ close enough to $t_0$\!;

\vspace{1ex}

\noindent $(ii)$\, $h^{0,\,1}_{\bar\partial}(t) = h^{0,\,1}_{\bar\partial}(t_0)$ and $h^{0,\,1}_{BC}(t) = h^{0,\,1}_{BC}(t_0)$ for all $t\in\Delta$ close enough to $t_0$.

\end{Cor}

 Unfortunately, the sGG property is not deformation closed (cf. Proposition \ref{not-closed}). However, another positive consequence of Theorem \ref{The:Betti-dbar} is the invariance of the sGG property under modifications.

\begin{Cor}\label{Cor:modification-invariance} Let $\mu\,:\,\widetilde{X}\to X$ be a holomorphic bimeromorphic map between compact complex manifolds $\widetilde{X}$ and $X$. The following equivalence holds\!\!:

$$\widetilde{X} \hspace{1ex}\mbox{is an sGG manifold} \hspace{3ex}\iff \hspace{3ex} X \hspace{1ex}\mbox{is an sGG manifold}.$$

\end{Cor}

 We would like to stress that although the sGG property of compact complex manifolds is a natural reinforcement (hopefully with multiple ramifications) of the strongly Gauduchon (sG) property, different reinforcements are possible. One of them is the following.

 We say that a Hermitian metric $\omega$ on a complex manifold $X$ of dimension $n$ is {\bf superstrong Gauduchon}\footnote{This term was coined by M. Verbitsky in a private communication with one of the authors who was simultaneously contemplating the same notion.} {\bf (super sG)} if $\partial\omega^{n-1}$ is $\partial\bar\partial$-exact, while $X$ is said to be a {\bf superstrong Gauduchon manifold (super sG manifold)} if it supports such a metric.

Any superstrong Gauduchon metric is trivially strongly Gauduchon and the two notions are equivalent if $X$ satisfies the $\partial\bar\partial$-lemma. It is also clear that any {\it balanced} metric $\omega$ (i.e. any $C^{\infty}$ positive definite $(1,\,1)$-form $\omega$ such that $\omega^{n-1}$ is $d$-closed, cf. [Gau77b] or [Mic83]) is superstrong Gauduchon. We sum up below the implication relations among these classes of compact complex manifolds $X$.

\noindent$\begin{array}{lllll} & & & &  \\
X\,\mbox{balanced manifold} & \implies & X\,\mbox{super sG manifold}  &  &   \\
 & \rotatebox{45}{$\implies$} & & \rotatebox{-45}{$\implies$} &  (\star\star) \\
X\,\partial\bar\partial\mbox{-manifold}  & \implies & X\,\mbox{sGG manifold} & \implies & X\,\mbox{sG manifold}\\
 & & & & \end{array}$

\noindent It will be proved in section \ref{section:sGG-other} that both implications on the last line in diagram $(\star\star)$ are strict, while the sGG and the super sG classes are unrelated (cf. Proposition \ref{unrelated-sGG-super-sG}) and so are the balanced and the sGG classes (cf. Proposition \ref{unrelated-sGG-balanced-Frolicher}). Thus, the class of sGG manifolds investigated in this work is new.

\vspace{3ex}

 $(2)$\, In connection with the second step of the approach to Conjecture \ref{Conj:classC_limit} outlined above, we prove in section \ref{section:cone-comparison} the following semi-continuity properties of the pseudo-effective and Gauduchon cones in families of sGG manifolds.

\begin{The}\label{The:semicontinuity-GE} Let $(X_t)_{t\in\Delta}$ be any holomorphic family of {\bf sGG} compact complex manifolds. Then ${\cal G}_{X_t}$ behaves lower-semicontinuously, while ${\cal E}_{X_t}$ behaves upper-semicontinuously w.r.t. the usual topology of $\Delta$ as $t\in\Delta$ varies.

\end{The}

 More precise statements will be given in Theorems \ref{The:lsc-G} and \ref{The:usc-E}. The main tool that we introduce in $\S.$\ref{subsection:fake} to prove Theorem \ref{The:semicontinuity-GE} in $\S.$\ref{subsection:s-c} is a pair of linear maps $(P,\,Q_{\omega})$ that we term {\it fake Hodge-Aeppli decomposition} of $H^{2n-2}_{DR}(X,\,\R)$ when $X$ is a compact sGG manifold of complex dimension $n$\!\!:

\vspace{1ex}

\hspace{10ex} $P\,:\,H^{2n-2}_{DR}(X,\,\R)\twoheadrightarrow H^{n-1,\,n-1}_A(X,\,\R)$

\vspace{1ex}

\noindent is a canonical surjection mimicking the projection onto $H^{n-1,\,n-1}_A(X,\,\R)$, while with any Hermitian metric $\omega$ on $X$ we associate a natural injection

\vspace{1ex}

\hspace{10ex} $Q_{\omega}\,:\, H^{n-1,\,n-1}_A(X,\,\R)\hookrightarrow H^{2n-2}_{DR}(X,\,\R)$

\vspace{1ex}

\noindent such that $P\circ Q_{\omega}$ is the identity map of $H^{n-1,\,n-1}_A(X,\,\R)$.

\vspace{2ex}

\noindent {\bf Notation.} On a given compact complex manifold $X$ of dimension $n$, for every $p,q=0,1, \dots , n$ we let $C^{\infty}_{p,\,q}(X,\,\C)$ stand for the space of $C^{\infty}$ forms of bidegree $(p,\,q)$, while $[\,\,]_{BC}, [\,\,]_A, [\,\,]_{\bar\partial}$ and $\{\,\,\}_{DR}$ will stand for Bott-Chern, Aeppli, Dolbeault and resp. De Rham cohomology classes. If a Hermitian metric $\omega$ has been fixed on $X$, we denote by $\Delta_A:=\partial\partial^{\star} + \bar\partial\bar\partial^{\star} + (\partial\bar\partial)^{\star}(\partial\bar\partial) + (\partial\bar\partial)(\partial\bar\partial)^{\star} + (\partial\bar\partial^{\star})(\partial\bar\partial^{\star})^{\star} + (\partial\bar\partial^{\star})^{\star}(\partial\bar\partial^{\star})$ the associated Aeppli Laplacian (which is a self-adjoint, elliptic operator of order $4$ introduced by Schweitzer in [Sch07] by analogy with the Bott-Chern Laplacian of Kodaira and Spencer) in which all the adjoints are computed w.r.t. the $L^2$ scalar product defined by $\omega$ on the space $C^{\infty}_{p,\,q}(X,\,\C)$ of $\C$-valued smooth $(p,\,q)$-forms on $X$. Thus $\Delta_A:C^{\infty}_{p,\,q}(X,\,\C)\to C^{\infty}_{p,\,q}(X,\,\C)$ and the Hodge isomorphism $H^{p,\,q}_A(X,\,\C)\simeq\ker\Delta_A$ holds, so every Aeppli cohomology class contains a unique $\Delta_A$-harmonic representative. (See e.g. [Pop13b, $\S.2$] for a review of these matters and further details.) On the other hand, $\Delta'':C^{\infty}_{p,\,q}(X,\,\C)\to C^{\infty}_{p,\,q}(X,\,\C)$ will denote the $\bar\partial$-Laplacian $\Delta''=\bar\partial\bar\partial^{\star} + \bar\partial^{\star}\bar\partial$.

\vspace{2ex}

\noindent {\bf Acknowledgments.} The authors are grateful to the referees for useful suggestions that contributed to the improvement of the presentation.

\section{Bott-Chern and Hodge numbers $h^{0,\,1}_{BC},$ $h^{0,\,1}_{\bar\partial}$}\label{section:BC-Hodge}

 In this section we prove a more precise version of Theorem \ref{The:BC-dbar}.

\begin{The}\label{The:S-map} Let $X$ be any compact complex manifold, $\mbox{dim}_{\C}X=n$.

\vspace{1ex}

\noindent $(i)$\, There is a well-defined canonical $\C$-linear map

 $$S\,:\,H^{n,\,n-1}_{\bar\partial}(X,\,\C)\longrightarrow H^{n,\,n-1}_A(X,\,\C), \hspace{2ex} S([\Gamma]_{\bar\partial}):=[\Gamma]_A.$$

\noindent Moreover, the map $S$ is {\bf surjective}, and we have an exact sequence

$$H^{n-1,\,n-1}_A(X,\,\C)\stackrel{T}{\longrightarrow}H^{n,\,n-1}_{\bar\partial}(X,\,\C)\stackrel{S}{\longrightarrow}H^{n,\,n-1}_A(X,\,\C)\longrightarrow 0,$$

\noindent i.e. $\mbox{Im}\,T = \ker S$, where $T$ is the map defined in (\ref{eqn:T-def}). In particular, $X$ is an sGG manifold if and only if $S$ is injective (i.e. if and only if $S$ is bijective).

\vspace{1ex}

\noindent $(ii)$\, There are well-defined canonical $\C$-linear maps and an exact sequence

$$0\longrightarrow H^{0,\,1}_{BC}(X,\,\C)\stackrel{S^{\star}}{\longrightarrow}H^{0,\,1}_{\bar\partial}(X,\,\C)\stackrel{T^{\star}}{\longrightarrow}H^{1,\,1}_{BC}(X,\,\C)$$

\noindent defined by $S^{\star}([u]_{BC}):=[u]_{\bar\partial}$ for any $d$-closed $(0,\,1)$-form $u$ and $T^{\star}([v]_{\bar\partial}):=[\partial v]_{BC}$ for any $\bar\partial$-closed $(0,\,1)$-form $v$. Thus $\mbox{Im}\,S^{\star} = \ker T^{\star}$.

Moreover, the maps $S^{\star}$ and $T^{\star}$ are dual to $S$ and respectively $T$. Thus $S^{\star}$ is {\bf injective}, hence $h^{0,\,1}_{BC}\leq h^{0,\,1}_{\bar\partial}$.

\vspace{1ex}

\noindent $(iii)$\, It follows that $X$ is an sGG manifold if and only if $S^{\star}$ is surjective (i.e. if and only if $S^{\star}$ is bijective) if and only if $h^{0,\,1}_{BC} = h^{0,\,1}_{\bar\partial}$.

\end{The}

\noindent {\it Proof.} By $S$ being well defined, we mean that $S([\Gamma]_{\bar\partial})$ (i.e. $[\Gamma]_A$) is meaningful and does not depend on the choice of representative $\Gamma$ of the class $[\Gamma]_{\bar\partial}$. The immediate verification of this fact is left to the reader.   

 It is clear that $\mbox{Im}\,T\subset\ker S$ since $\mbox{Im}\,\partial\subset\mbox{Im}\,\partial + \mbox{Im}\,\bar\partial$. To show the reverse inclusion, let $[\Gamma]_{\bar\partial}\in H^{n,\,n-1}_{\bar\partial}(X,\,\C)$ such that $S([\Gamma]_{\bar\partial})=0$. Then there are forms $\Omega, \Lambda$ of respective bidegrees $(n-1,\,n-1)$ and $(n,\,n-2)$ such that $\Gamma = \partial\Omega + \bar\partial\Lambda$, i.e. $\Gamma - \bar\partial\Lambda = \partial\Omega$. Hence $\partial\bar\partial\Omega = 0$, $[\Gamma]_{\bar\partial} = [\Gamma - \bar\partial\Lambda]_{\bar\partial}$ and $T([\Omega]_A) = [\partial\Omega]_{\bar\partial} = [\Gamma]_{\bar\partial}$. Thus $[\Gamma]_{\bar\partial}\in\mbox{Im}\,T$. This proves the identity $\mbox{Im}\,T = \ker S$.

 The surjectivity of $S$ will follow from the injectivity of its dual map $S^{\star}$ that will be proved below.

 The well-definedness of $S^{\star}$ and $T^{\star}$ are proved in a similar way. The identity $\mbox{Im}\,S^{\star} = \ker T^{\star}$ follows by duality from $\mbox{Im}\,T = \ker S$ or directly in the following way. Let $[u]_{BC}\in H^{0,\,1}_{BC}(X,\,\C)$, i.e. $u$ is a $d$-closed $(0,\,1)$-form. Then $\partial u=0$ and $\bar\partial u=0$, hence $T^{\star}(S^{\star}[u]_{BC}) = T^{\star}([u]_{\bar\partial}) = [\partial u]_{BC} = 0$. Thus $\mbox{Im}\,S^{\star}\subset\ker T^{\star}$. To show the reverse inclusion, let $[v]_{\bar\partial}\in\ker T^{\star}$, i.e. $v$ is a $\bar\partial$-closed $(0,\,1)$-form such that $\partial v = \partial\bar\partial f$ for some function $f$. Then $\partial(v - \bar\partial f) =0$, hence $d(v - \bar\partial f) = 0$ and $[v]_{\bar\partial} = [v - \bar\partial f]_{\bar\partial} = S^{\star}([v - \bar\partial f]_{BC})$, so $[v]_{\bar\partial}\in\mbox{Im}\,S^{\star}$.

 Let us now show that $S^{\star}$ is injective. Let $[u]_{BC}\in\ker S^{\star}$, i.e. $u$ is a $d$-closed $(0,\,1)$-form such that $u = \bar\partial f$ for some function $f$. Since $du=0$, we also have $\partial u=0$, hence $\partial\bar\partial f = 0$ on $X$. Because $X$ is compact, the function $f$ must be constant, hence $u = \bar\partial f = 0$. In particular, $[u]_{BC} = 0$.

Let us now check that the maps $T$ and $T^{\star}$ are dual to each other under the duality (\ref{eqn:BC-A-duality}) and under the Serre duality

\vspace{1ex}

\hspace{3ex} $H^{0,\,1}_{\bar\partial}(X,\,\C)\times H^{n,\,n-1}_{\bar\partial}(X,\,\C)\longrightarrow\C, \hspace{3ex} ([v]_{\bar\partial},\,[\Gamma]_{\bar\partial})\mapsto\int\limits_Xv\wedge\Gamma,$

\vspace{1ex}

\noindent the latter being defined for every $\bar\partial$-closed forms $v$ and $\Gamma$ of respective bidegrees $(0,\,1)$ and $(n,\,n-1)$. We have to check that for every $[v]_{\bar\partial}\in H^{0,\,1}_{\bar\partial}(X,\,\C)\simeq (H^{n,\,n-1}_{\bar\partial}(X,\,\C))^{\star}$, if we denote by

\vspace{1ex}

\hspace{6ex} $\sigma_v\,:\,H^{n,\,n-1}_{\bar\partial}(X,\,\C)\longrightarrow\C$

\vspace{1ex}

\noindent the linear map induced by $[v]_{\bar\partial}$ under duality, then the linear map

\vspace{1ex}

\hspace{6ex} $\tau_{\partial v}\,:\,H^{n-1,\,n-1}_A(X,\,\C)\longrightarrow\C$

\vspace{1ex}

\noindent induced by $T^{\star}([v]_{\bar\partial}) = [\partial v]_{BC}\in H^{1,\,1}_{BC}(X,\,\C)\simeq (H^{n-1,\,n-1}_A(X,\,\C))^{\star}$ under duality is $\sigma_v\circ T$.
This is indeed the case since, for every $[\Omega]_A\in H^{n-1,\,n-1}_A(X,\,\C)$, we have

\vspace{1ex}

\hspace{1ex} $\displaystyle(\sigma_v\circ T)([\Omega]_A) = \sigma_v([\partial\Omega]_{\bar\partial}) = \int\limits_Xv\wedge\partial\Omega = \int\limits_X\partial v\wedge\Omega = \tau_{\partial v}([\Omega]_A),$

\vspace{1ex}

\noindent having used the Stokes formula $\int_X\partial(v\wedge\Omega) = 0$ and $\partial(v\wedge\Omega) = \partial v\wedge\Omega - v\wedge\partial\Omega$.

 We can now check the equivalence\!\!:

\vspace{1ex}

\hspace{6ex} $T$ vanishes identically $\iff$ $T^{\star}$ vanishes identically.

\vspace{1ex}

\noindent Indeed, $T$ vanishes identically if and only if for every $[\Omega]_A\in H^{n-1,\,n-1}_A(X,\,\C)$ and every $[v]_{\bar\partial}\in H^{0,\,1}_{\bar\partial}(X,\,\C)$ we have $\int_X\partial\Omega\wedge v = 0$. Since $\int_X\partial\Omega\wedge v = \int_X\Omega\wedge\partial v$ by the Stokes formula, this is equivalent to the map $\tau_{\partial v}\,:\,H^{n-1,\,n-1}_A(X,\,\C)\rightarrow\C$ vanishing identically, i.e. to $[\partial v]_{BC} = 0$, for every $[v]_{\bar\partial}\in H^{0,\,1}_{\bar\partial}(X,\,\C)$. Since $[\partial v]_{BC} = T^{\star}([v]_{\bar\partial})$, this is still equivalent to the map $T^{\star}$ vanishing identically.

 Thus, if we put the various bits together, we get the equivalences\!\!:

\vspace{1ex}

\noindent $X\hspace{1ex}\mbox{is an sGG manifold} \iff \mbox{Im}\,T = 0 \iff \ker S = 0 \iff S \hspace{1ex} \mbox{is injective}$

\vspace{1ex}

\hspace{19ex} $\iff \ker T^{\star} = H^{0,\,1}_{\bar\partial}(X,\,\C) \iff \mbox{Im}\,S^{\star} = H^{0,\,1}_{\bar\partial}(X,\,\C)$

\vspace{1ex}

\hspace{19ex} $\iff S^{\star} \hspace{1ex} \mbox{is surjective}.$

\vspace{1ex}

 \noindent It can be checked that the maps $S$ and $S^{\star}$ are dual to each other in the same way as the duality between $T$ and $T^{\star}$ has been checked. \hfill $\Box$

\vspace{3ex}

\noindent {\it Proof of Corollary \ref{Cor:Iwasawa-sGG}.} Reading the dimension tables for the Hodge and Bott-Chern numbers given in [Nak75, p.96] and resp. [Ang11, Theorem 5.1], we gather that

$$h^{0,\,1}_{BC} = h^{0,\,1}_{\bar\partial} = 2$$

\noindent for the Iwasawa manifold and all its small deformations. Thus, the conclusion follows from Theorem \ref{The:BC-dbar}.  \hfill $\Box$

\begin{Rem}\label{solvmanifolds}
{\rm 
In the context of solvmanifolds some examples of sGG manifolds can be obtained.
For instance, for the completely-solvable Nakamura manifold, studied first by Nakamura in [Nak75], it is shown by Angella 
and Kasuya that the corresponding Lie group $G$ admits lattices $\Gamma$ (see cases (i)--(iii) in [AK12, Example 2.17])
for which the Bott-Chern cohomology of 
the compact solvmanifolds $G/\Gamma$ can be determined. For the lattices $\Gamma$ in cases (ii) and (iii) 
the solvmanifolds satisfy $h^{0,\,1}_{BC}(G/\Gamma)=1=h^{0,\,1}_{\bar\partial}(G/\Gamma)$ (see [AK12, Table 6]),
so by Theorem \ref{The:BC-dbar} they are sGG.
Note that in [AK12, Remark 2.19] it is proved that $G/\Gamma$ is not a $\partial\bar\partial$-manifold only for $\Gamma$ in case (ii).
}
\end{Rem}

\section{Betti and Hodge numbers $b_1$, $h^{0,\,1}_{\bar\partial}$}\label{section:Betti-Hodge}

  In this section we prove a more precise version of Theorem \ref{The:Betti-dbar} from which the latter follows as a corollary. For any form $\alpha$, we denote by $\alpha^{p,\,q}$ its component of bidegree $(p,\,q)$.

\begin{The}\label{The:main} Let $X$ be any compact complex manifold, $\mbox{dim}_{\C}X=n$.

\vspace{1ex}

\noindent $(i)$\, There is a well-defined canonical $\C$-linear map

\begin{eqnarray}\nonumber F\,:\,H^1_{DR}(X,\,\C) & \longrightarrow & H^{0,\,1}_{\bar\partial}(X,\,\C)\oplus\overline{H^{0,\,1}_{\bar\partial}(X,\,\C)},\\
\nonumber  F(\{\alpha\}_{DR}) & := & ([\alpha^{0,\,1}]_{\bar\partial},\,\overline{[\overline{\alpha^{1,\,0}}]_{\bar\partial}}).\end{eqnarray}

\noindent Moreover, the map $F$ is {\bf injective}. Consequently, the following inequality holds on any compact complex manifold\!\!:

$$b_1\leq 2h^{0,\,1}_{\bar\partial}.$$

\vspace{1ex}

\noindent $(ii)$\, There is a well-defined canonical $\C$-linear map\!\!:

\begin{eqnarray}\nonumber F^{\star}:\,H^{n,\,n-1}_{\bar\partial}(X,\,\C)\oplus\overline{H^{n,\,n-1}_{\bar\partial}(X,\,\C)} & \longrightarrow & H^{2n-1}_{DR}(X,\,\C),\\
\nonumber  F^{\star}([\beta]_{\bar\partial},\,\overline{[\gamma]}_{\bar\partial}) & := & \{\beta + \bar\gamma\}_{DR}.\end{eqnarray}

\noindent Moreover, the map $F^{\star}$ is dual to the map $F$. Hence $F^{\star}$ is {\bf surjective}.

\vspace{1ex}

\noindent $(iii)$\, The following equivalence holds\!\!:

$$X\hspace{1ex}\mbox{is an sGG manifold} \iff F^{\star}\hspace{1ex}\mbox{is injective}.$$

\noindent Since $F$ is always injective by $(i)$, this means that $X$ is an sGG manifold if and only if the linear map $F$ is bijective. In other words, the following equivalence holds\!\!:

$$X\hspace{1ex}\mbox{is an sGG manifold} \iff b_1 = 2\,h^{0,\,1}_{\bar\partial}.$$

\end{The}

\noindent {\it Proof.} $(i)$\, For any $1$-form $\alpha$, the condition $d\alpha=0$ is equivalent to

$$\bar\partial\alpha^{0,\,1} = 0, \hspace{2ex} \partial\alpha^{1,\,0}=0\hspace{1ex} (\Leftrightarrow \bar\partial\overline{\alpha^{1,\,0}}=0), \hspace{2ex} \partial\alpha^{0,\,1} + \bar\partial\alpha^{1,\,0}=0.$$

\noindent Thus, if $d\alpha=0$, $\alpha^{0,\,1}$ and $\overline{\alpha^{1,\,0}}$ define Dolbeault cohomology classes of type $(0,\,1)$. To show that the map $F$ is independent of the choice of representative in a given De Rham class $\{\alpha\}_{DR}$, let $\alpha$ be any $d$-exact $1$-form on $X$. Then, there exists a function $f$ on $X$ such that $\alpha = df = \partial f + \bar\partial f$. Hence $\alpha^{0,\,1}=\bar\partial f$ and $\overline{\alpha^{1,\,0}}=\bar\partial\bar{f}$, so $[\alpha^{0,\,1}]_{\bar\partial} = [\overline{\alpha^{1,\,0}}]_{\bar\partial} = 0$. This proves the well-definedness of the map $F$.

 To prove that $F$ is injective, let $\alpha$ be a $d$-closed $1$-form such that $F(\{\alpha\}_{DR}) = 0$, i.e. $\alpha^{0,\,1} = \bar\partial f$ and $\overline{\alpha^{1,\,0}} = \bar\partial g$ (i.e. $\alpha^{1,\,0} = \partial\bar{g}$) for some functions $f, g$ on $X$. Then

$$0 = \partial\alpha^{0,\,1} + \bar\partial\alpha^{1,\,0} = \partial\bar\partial(f-\bar{g})  \hspace{2ex} \mbox{on}\hspace{1ex} X,$$

\noindent where the first identity follows from $\partial\alpha^{0,\,1} + \bar\partial\alpha^{1,\,0}$ being the component of bidegree $(1,\,1)$ of $d\alpha = 0$. Since $X$ is compact, $f-\bar{g}$ must be constant on $X$, hence $\partial\bar{g} = \partial f$, so we get

$$\alpha = \alpha^{1,\,0} + \alpha^{0,\,1}  = \partial f + \bar\partial f = df.$$

\noindent Thus $\{\alpha\}_{DR} = 0$. Consequently, $F$ is injective.

\vspace{1ex}

\noindent $(ii)$\, For any $\bar\partial$-closed $(n,\,n-1)$-forms $\beta, \gamma$, we have $\partial\beta = \partial\gamma = 0$ for bidegree reasons, hence $d\beta = d\gamma = 0$, so $\beta + \bar{\gamma}$ is $d$-closed and therefore it defines a De Rham class. To show that $F^{\star}$ is independent of the choice of representatives of the classes $[\beta]_{\bar\partial}, [\gamma]_{\bar\partial}$, suppose that $\beta = \bar\partial u$ and $\gamma = \bar\partial v$ for some $(n,\,n-2)$-forms $u, v$. Since $\partial u = \partial v = 0$ for bidegree reasons, we see that $\beta = du$ and $\gamma = dv$, hence $\beta + \bar{\gamma} = d(u + \bar{v})$, so $\{\beta + \bar{\gamma}\}_{DR} = 0$. We conclude that $F^{\star}$ is well defined.

 We now prove that the maps $F$ and $F^{\star}$ are dual to each other. Under the Serre duality $H^{n,\,n-1}_{\bar\partial}(X,\,\C)\simeq(H^{0,\,1}_{\bar\partial}(X,\,\C))^{\star}$, every pair $([\beta]_{\bar\partial},\,\overline{[\gamma]_{\bar\partial}})\in H^{n,\,n-1}_{\bar\partial}(X,\,\C)\oplus\overline{H^{n,\,n-1}_{\bar\partial}(X,\,\C)}$ can be identified with the pair $(u,\,\bar{v})$ in which $u, v\,:\,H^{0,\,1}_{\bar\partial}(X,\,\C)\rightarrow\C$ are the $\C$-linear maps acting as

$$u([\alpha^{0,\,1}]_{\bar\partial}) = \int\limits_X\beta\wedge\alpha^{0,\,1} \hspace{2ex} \mbox{and} \hspace{2ex} v([\alpha^{0,\,1}]_{\bar\partial}) = \int\limits_X\gamma\wedge\alpha^{0,\,1}$$

\noindent for every class $[\alpha^{0,\,1}]_{\bar\partial}\in H^{0,\,1}_{\bar\partial}(X,\,\C)$ and $\bar{v}\,:\,\overline{H^{0,\,1}_{\bar\partial}(X,\,\C)}\to\C$ is the $\C$-linear map defined by $\bar{v}(\overline{[\alpha^{0,\,1}]_{\bar\partial}}):=\overline{v([\alpha^{0,\,1}]_{\bar\partial})}$. Proving the duality between $F$ and $F^{\star}$ amounts to proving that

\begin{equation}\label{eqn:dualityFstar}\sigma_{\beta + \bar\gamma} = (u + \bar{v})\circ F\end{equation}

\noindent for any $[\beta]_{\bar\partial}, [\gamma]_{\bar\partial}\in H^{n,\,n-1}_{\bar\partial}(X,\,\C)$, where $\sigma_{\beta + \bar\gamma}\,:\,H^1_{DR}(X,\,\C)\rightarrow\C$ is the $\C$-linear map representing the De Rham class

$$F^{\star}([\beta]_{\bar\partial},\,\overline{[\gamma]_{\bar\partial}}) = \{\beta + \bar{\gamma}\}_{DR}\in H^{2n-1}_{DR}(X,\,\C)\simeq(H^1_{DR}(X,\,\C))^{\star}$$

\noindent under Poincar\'e duality. By definition, this means that for every $\{\alpha\}_{DR}\in H^1_{DR}(X,\,\C)$ we have

$$\sigma_{\beta + \bar\gamma}(\{\alpha\}_{DR}) = \int\limits_X(\beta + \bar\gamma)\wedge\alpha = \int\limits_X\beta\wedge\alpha^{0,\,1} + \int\limits_X\bar\gamma\wedge\alpha^{1,\,0} = \bigg((u + \bar{v})\circ F\bigg)(\{\alpha\}_{DR}).$$

\noindent This proves (\ref{eqn:dualityFstar}). We conclude that $F$ and $F^{\star}$ are dual to each other.

\vspace{1ex}

\noindent $(iii)$\, Let us first prove the implication ``$\implies$''. Suppose that $X$ is an sGG manifold. Let $[\beta]_{\bar\partial}, [\gamma]_{\bar\partial}\in H^{n,\,n-1}_{\bar\partial}(X,\,\C)$ such that $\beta + \bar{\gamma} = d(\Omega^{n,\,n-2} + \Omega^{n-1,\,n-1} + \Omega^{n-2,\,n})$ for some forms $\Omega^{p,\,q}$ of the specified bidegrees. This amounts to having

$$\beta = \bar\partial\Omega^{n,\,n-2} + \partial\Omega^{n-1,\,n-1} \hspace{2ex}\mbox{and}\hspace{2ex} \bar{\gamma} = \bar\partial\Omega^{n-1,\,n-1} + \partial\Omega^{n-2,\,n},$$

\noindent therefore to having

$$\beta - \bar\partial\Omega^{n,\,n-2} = \partial\Omega^{n-1,\,n-1} \hspace{2ex}\mbox{and}\hspace{2ex} \gamma - \bar\partial\overline{\Omega^{n-2,\,n}} = \partial\overline{\Omega^{n-1,\,n-1}}.$$

\noindent Now, $\beta - \bar\partial\Omega^{n,\,n-2}\in[\beta]_{\bar\partial}$ and $\gamma - \bar\partial\overline{\Omega^{n-2,\,n}}\in[\gamma]_{\bar\partial}$. On the other hand, $\partial\Omega^{n-1,\,n-1}$ and $\partial\overline{\Omega^{n-1,\,n-1}}$ are $d$-closed and $\partial$-exact $(n,\,n-1)$-forms on $X$, so the sGG assumption on $X$ implies (thanks to $(iii)$ of Lemma \ref{Lem:sGG-equiv-descriptions1}) that they are both $\bar\partial$-exact, i.e. $[\beta]_{\bar\partial} = [\gamma]_{\bar\partial} = 0$ in $H^{n,\,n-1}_{\bar\partial}(X,\,\C)$. We conclude that $F^{\star}$ is injective if $X$ is sGG.

 We now prove the reverse implication ``$\Longleftarrow$''. Suppose that $F^{\star}$ is injective. We will show that every Gauduchon metric on $X$ is actually strongly Gauduchon. This will imply that $X$ is an sGG manifold thanks to $(iv)$ of Lemma \ref{Lem:sGG-equiv-descriptions1}.

 Let $\omega$ be any Gauduchon metric on $X$. Then $\partial\omega^{n-1}\in\ker\bar\partial$, so we have a Dolbeault class $[\partial\omega^{n-1}]_{\bar\partial}\in H^{n,\,n-1}_{\bar\partial}(X,\,\C)$. Now, $\partial\omega^{n-1} + \overline{\partial\omega^{n-1}} = d\omega^{n-1}$ and, moreover,

$$F^{\star}([\partial\omega^{n-1}]_{\bar\partial},\,\overline{[\partial\omega^{n-1}]_{\bar\partial}}) = \{\partial\omega^{n-1} + \overline{\partial\omega^{n-1}}\}_{DR} = \{d\omega^{n-1}\}_{DR} = 0.$$

\noindent Since $F^{\star}$ is supposed injective, we infer that $[\partial\omega^{n-1}]_{\bar\partial} = 0$, i.e. $\omega$ is strongly Gauduchon.  \hfill $\Box$

\vspace{3ex}

\noindent {\it Proof of Corollary \ref{Cor:openness}.} This is an immediate consequence of Theorem \ref{The:Betti-dbar} if we use the (local) invariance of the Betti numbers of the fibres in a $C^{\infty}$ family of compact complex manifolds and the upper-semicontinuity of the Hodge numbers $h^{p,\,q}(t)$ as $t$ varies in $\Delta$ (cf. [KS60, Theorem 4]). Indeed, if $X_{t_0}$ is an sGG manifold, we have\!\!:

$$b_1 = 2h^{0,\,1}_{\bar\partial}(t_0)\geq 2h^{0,\,1}_{\bar\partial}(t)\geq b_1 \hspace{2ex}\mbox{for all}\hspace{1ex}t\hspace{1ex}\mbox{sufficiently close to}\hspace{1ex}t_0.$$

\noindent Thus, we must have equalities $b_1 = 2h^{0,\,1}_{\bar\partial}(t_0)= 2h^{0,\,1}_{\bar\partial}(t)$ for all $t$ close to $t_0$. In particular, by Theorem \ref{The:Betti-dbar}, $X_t$ must be an sGG manifold for all $t$ close to $t_0$. Then Theorem \ref{The:BC-dbar} implies $h^{0,\,1}_{BC}(t_0)=h^{0,\,1}_{BC}(t)$ for $t$ close to $t_0$. \hfill $\Box$

\vspace{3ex}

\noindent {\it Proof of Corollary \ref{Cor:modification-invariance}.} Thanks to the characterisation of sGG manifolds given in Theorem \ref{The:Betti-dbar}, it suffices to ensure the invariance of $b_1$ and of $h^{0,\,1}_{\bar\partial}$ under modifications, both of which are standard. Indeed, the fundamental group is known to be a bimeromorphic invariant of complex manifolds, hence so is its abelianisation $H_1$, so also $b_1$.

 We recall for the reader's convenience the well-known argument showing the modification invariance of every $h^{0,\,k}_{\bar\partial}$ for any compact complex manifold (not necessarily satisfying the Hodge symmetry)\footnote{One of the authors is very grateful to Professors A. Fujiki and F. Campana for pointing out to him this simple argument.}. This invariance follows from the combination of two things. The first thing is the following standard fact (cf. e.g. [Har77]) giving the vanishing of the higher direct image sheaves of the structural sheaf under modifications\!\!:

\vspace{1ex}

{\it Let $f\,:\,X\to Y$ be a {\bf bimeromorphic} morphism between (smooth) compact complex manifolds. Then\!\!:

\vspace{1ex}

\noindent $(i)$\, $f_{\star}{\cal O}_X = {\cal O}_Y$\!\!;

\vspace{1ex}

\noindent $(ii)$\,$R^if_{\star}{\cal O}_X = 0$ for all $i>0$.}

\vspace{1ex}

\noindent The second thing is the Leray spectral sequence associated with $f$ and ${\cal O}_X$. Recall that this is the spectral sequence starting at $E_2^{p,\,q}:=H^p(Y,\,R^qf_{\star}{\cal O}_X)$ and converging to $H^{p+q}(X,\,{\cal O}_X)$. The shape of the direct image sheaves of ${\cal O}_X$ under $f$ implies at once that

$$E_2^{p,\,0} = H^p(Y,\,{\cal O}_Y)\simeq H^{0,\,p}(Y,\,\C) \hspace{2ex} \mbox{and} \hspace{2ex} E_2^{p,\,q} = 0, \hspace{1ex} q\geq 1.$$

\noindent It follows that the Leray spectral sequence degenerates at $E_2$ and we have

$$H^{0,\,k}(X,\,\C) = H^k(X,\,{\cal O}_X)\simeq\bigoplus_{p+q=k}E_2^{p,\,q} = E_2^{k,\,0}\simeq H^{0,\,k}(Y,\,\C)   \hspace{3ex} \mbox{for all}\hspace{1ex} k.  $$

\hfill $\Box$

\section{The interior of the pseudo-effective cone}\label{section:interior-psef}

Starting from a handful of trivial observations, we exhibit in this section a few examples of compact complex manifolds which are not in the class ${\cal C}$ but whose pseudo-effective cone has non-empty interior.

\begin{Prop}\label{Prop:nonemptyE} $(I)$\, Let $X$ be a compact complex {\bf surface}. Then\!\!:

\vspace{1ex}

\noindent $(i)$\,there exists a non-zero $d$-closed $(1,\,1)$-current $T\geq 0$ on $X$\!;

\vspace{1ex}

\noindent $(ii)$\,$h^{1,\,1}_{BC}(X,\,\C)\geq 1$\!;

\vspace{1ex}

\noindent $(iii)$\, if $h^{1,\,1}_{BC}(X,\,\C) = 1$, then $\mathring{{\cal E}}_X\neq\emptyset$.

\vspace{1ex}

 $(II)$\, Let $X$ be a compact complex manifold of any dimension.

\vspace{1ex}

\noindent $(i)$\, If $h^{1,\,1}_{BC}(X,\,\C) = 1$ and if there exists a non-zero $d$-closed $(1,\,1)$-current $T\geq 0$ on $X$, then $\mathring{{\cal E}}_X\neq\emptyset$.

\vspace{1ex}

\noindent $(ii)$\, If $b_2(X)=0$, then $X$ is not in the class ${\cal C}$.

\end{Prop}

\noindent {\it Proof.} $(I)(i)$\, Suppose that such a current did not exist. Then by [Pop13b, Proposition 5.4], there would exist a degenerate balanced structure $\omega$ on $X$, i.e. a smooth $(1,\,1)$-form $\omega>0$ such that $\omega^{n-1}$ is $d$-exact. Since $n-1 = 1$ on a surface, $\omega$ would be, in particular, a K\"ahler metric, contradicting the supposed non-existence of a non-zero $d$-closed positive $(1,\,1)$-current.

 $(ii)$\, Let $T\geq 0$ be a non-zero $d$-closed $(1,\,1)$-current on $X$ (which exists by $(i)$). Then $[T]_{BC}\in H^{1,\,1}_{BC}(X,\,\C)$ cannot be the zero Bott-Chern class since, otherwise, $T=i\partial\bar\partial\varphi\geq 0$ on $X$ for some $L^1_{loc}$ function $\varphi$, so $\varphi$ would be a global psh function on the compact manifold $X$, hence $\varphi$ would be constant and $T=i\partial\bar\partial\varphi=0$, a contradiction. 

$(iii)$\, For any non-zero $d$-closed $(1,\,1)$-current $T\geq 0$ on $X$, we have $0\neq [T]_{BC}\in{\cal E}_X\subset H^{1,\,1}_{BC}(X,\,\R)$. Since ${\cal E}_X$ is a convex cone, it must contain the whole ray $\R^{+}\!\!\cdot\![T]_{BC}$, so it has non-empty interior in the ambient $1$-dimensional real vector space.

\vspace{1ex}

$(II)(i)$\, The proof of this statement is identical to that of $(I)(iii)$. It has been necessary to suppose the existence of a non-zero $d$-closed positive $(1,\,1)$-current since, unlike compact complex surfaces, arbitrary compact complex manifolds of dimension $\geq 3$ need not possess such a current (see e.g. [Pop13b]).

$(II)(ii)$\, Suppose that $X$ is in the class ${\cal C}$. Then, on the one hand, $X$ is an sG manifold, while on the other hand, by [DP04], there exists a K\"ahler current $T$ on $X$. In particular, $T$ defines a De Rham class $\{T\}_{DR}\in H^2_{DR}(X,\,\R)=0$, hence $T$ is a $d$-exact (non-zero positive) $(1,\,1)$-current. The existence of a such a current is equivalent, thanks to [Pop13a, Proposition 4.3], to the manifold $X$ not being an sG manifold, a contradiction.  \hfill $\Box$

\vspace{2ex}

 We now notice a few examples showing that the property $\mathring{{\cal E}}_X\neq\emptyset$ does not imply that $X$ is a class ${\cal C}$ manifold.

\begin{Prop}\label{Prop:nonemptyE-examples} Let $X$ be either a {\bf Hopf surface}, or an {\bf Inoue $S_M$ surface}, or {\bf an Inoue $S_{\pm}$ surface}, or a {\bf secondary Kodaira surface}.

 Then $X$ is not in the class ${\cal C}$ but $\mathring{{\cal E}}_X\neq\emptyset$.

\end{Prop}

\noindent {\it Proof.} All the surfaces of the above types are non-K\"ahler, hence not in the class ${\cal C}$ (since the K\"ahler class coincides with the class ${\cal C}$ in the case of surfaces). Now, thanks to [ADT14, Theorem 2.2, Tables 1 and 2, p.8-9]\footnote{We are very grateful to F. Campana for pointing out this reference to us.}, $h^{1,\,1}_{BC}(X,\,\C)=1$ for each of these surfaces. We get $\mathring{{\cal E}}_X\neq\emptyset$ from part $(I)(iii)$ of our Proposition \ref{Prop:nonemptyE}.  \hfill $\Box$

\section{The cones ${\cal G}_X$ and ${\cal E}_X$ under deformations}\label{section:cone-comparison}

 Throughout this section, as in the rest of the paper, for any differential form $\Omega$ of any degree $k$ and for any $(p,\,q)$ such that $p+q=k$, we denote by $\Omega^{p,\,q}$ the component of $\Omega$ of bidegree $(p,\,q)$. Thus $\Omega = \sum\limits_{p+q=k}\Omega^{p,\,q}$.

 It will be seen in $\S.$\ref{subsection:s-c} that the discussion of the variation of the cones ${\cal G}_X$ and ${\cal E}_X$ under deformations of the complex structure of a compact sGG manifold $X$ would be greatly simplified if the Bott-Chern number $h^{1,\,1}_{BC}(X)$ were locally deformation constant. Unfortunately, this is not the case as Proposition~\ref{Prop:h11-jumping} will show, rendering indispensable the introduction of some technical aspects in $\S.$\ref{subsection:s-c}.

\subsection{Fake Hodge-Aeppli decomposition of $H^{2n-2}_{DR}(X,\,\R)$}\label{subsection:fake}

 If our $n$-dimensional compact complex manifold $X$ were supposed to be a $\partial\bar\partial$-manifold, there would exist a canonical isomorphism $H^{2n-2}_{DR}(X,\,\C)\simeq H^{n,\,n-2}_A(X,\,\C)\oplus H^{n-1,\,n-1}_A(X,\,\C)\oplus H^{n-2,\,n}_A(X,\,\C)$ (cf. [Pop13b] where this splitting was called a Hodge-Aeppli decomposition), hence in particular a canonical surjection $H^{2n-2}_{DR}(X,\,\C)\twoheadrightarrow H^{n-1,\,n-1}_A(X,\,\C)$ and a canonical injection $H^{n-1,\,n-1}_A(X,\,\C)\hookrightarrow H^{2n-2}_{DR}(X,\,\C)$ which is a section of the surjection. However, under the weaker sGG assumption on $X$, a complete Hodge-Aeppli decomposition in degree $2n-2$ need not exist, but we will show that a weaker substitute thereof (that will prove sufficient for our purposes later on) exists\!: a canonical surjection and a non-canonical but naturally-associated-with-any-given-metric injection as above exist if we restrict attention to the real cohomologies.

 We start by noticing the existence of the canonical surjection.

\begin{Prop}\label{Prop:P-def-surj} Let $X$ be an arbitrary compact complex manifold, $\mbox{dim}_{\C}X=n$. The following canonical linear map

\begin{equation}\label{eqn:P-def-surj} P\,:\,H^{2n-2}_{DR}(X,\,\R)\to H^{n-1,\,n-1}_A(X,\,\R), \hspace{2ex} \{\Omega\}_{DR}\mapsto [\Omega^{n-1,\,n-1}]_A,\end{equation}

\noindent is well defined. Furthermore, if $X$ is an sGG manifold, $P$ is surjective.

\end{Prop}

\noindent {\it Proof.} Let $\Omega = \Omega^{n,\,n-2} + \Omega^{n-1,\,n-1} + \Omega^{n-2,\,n}$ be any $d$-closed (not necessarily real) $C^{\infty}$ form of degree $2n-2$. We have\!\!:

\begin{eqnarray}\label{eqn:d-closedness-cond}\nonumber d\Omega = 0 & \iff & \partial\Omega^{n-1,\,n-1} + \bar\partial\Omega^{n,\,n-2} = 0 \hspace{1ex}\mbox{and}\hspace{1ex} \partial\Omega^{n-2,\,n} + \bar\partial\Omega^{n-1,\,n-1} = 0\\
   & \implies & \partial\bar\partial\Omega^{n-1,\,n-1} =0.\end{eqnarray}

\noindent The last identity shows that $\Omega^{n-1,\,n-1}$ defines indeed an Aeppli cohomology class of bidegree $(n-1,\,n-1)$. To show well-definedness for $P$, we still have to show that the definition is independent of the choice of representative of the De Rham class $\{\Omega\}_{DR}$. Let $\Omega_1, \Omega_2$ represent the same De Rham class, i.e. $\Omega:=\Omega_1 - \Omega_2$ is $d$-exact. Then there exists a $(2n-3)$-form $\Gamma$ such that $\Omega = d(\Gamma^{n,\,n-3} + \Gamma^{n-1,\,n-2} + \Gamma^{n-2,\,n-1} + \Gamma^{n-3,\,n})$. Thus, the $d$-exactness of a $(2n-2)$-form $\Omega$ is equivalent to the existence of $\Gamma\in C^{\infty}_{2n-3}(X,\,\C)$ such that

\begin{eqnarray}\label{eqn:d-exactness-cond}\nonumber \Omega^{n,\,n-2} & \stackrel{(i)}{=} & \partial\Gamma^{n-1,\,n-2} + \bar\partial\Gamma^{n,\,n-3}\\
\nonumber \Omega^{n-1,\,n-1} & \stackrel{(ii)}{=} & \partial\Gamma^{n-2,\,n-1} + \bar\partial\Gamma^{n-1,\,n-2}\\
 \Omega^{n-2,\,n} & \stackrel{(iii)}{=} & \partial\Gamma^{n-3,\,n} + \bar\partial\Gamma^{n-2,\,n-1}.\end{eqnarray}

\noindent Identity $(ii)$ above means that $[\Omega^{n-1,\,n-1}]_A = 0$, i.e. $[\Omega^{n-1,\,n-1}_1]_A = [\Omega^{n-1,\,n-1}_2]_A$.

 Let us now suppose that $X$ is an sGG manifold. Pick any class $[\Omega^{n-1,\,n-1}]_A\in H^{n-1,\,n-1}_A(X,\,\R)$ with $\Omega^{n-1,\,n-1}$ {\it real}. Then $\partial\bar\partial\Omega^{n-1,\,n-1}=0$, hence $d(\partial\Omega^{n-1,\,n-1})=0$. Since $\partial\Omega^{n-1,\,n-1}$ is a $\partial$-exact $d$-closed $(n,\,n-1)$-form, the sGG assumption on $X$ implies that $\partial\Omega^{n-1,\,n-1}$ is $\bar\partial$-exact (see Lemma \ref{Lem:sGG-equiv-descriptions1}). Thus, there exists an $(n,\,n-2)$-form $\Omega^{n,\,n-2}$ such that

 $$\partial\Omega^{n-1,\,n-1} = - \bar\partial\Omega^{n,\,n-2}.$$

\noindent Hence, since $\Omega^{n-1,\,n-1}$ is real, $\bar\partial\Omega^{n-1,\,n-1} = - \partial\overline{\Omega^{n,\,n-2}}$. Therefore, the $(2n-2)$-form $\Omega:=\Omega^{n,\,n-2} + \Omega^{n-1,\,n-1} + \overline{\Omega^{n,\,n-2}}$ is real and $d\Omega = 0$ (cf. (\ref{eqn:d-closedness-cond})). It is clear that $P(\{\Omega\}_{DR}) = [\Omega^{n-1,\,n-1}]_A$. This proves that $P$ is surjective.  \hfill $\Box$

\begin{Cor}\label{Cor:Pstar-inj} If $X$ is an sGG compact complex manifold with $\mbox{dim}_{\C}X=n$, the dual map of $P$\!\!:

\begin{equation}\label{eqn:Pstar-inj} P^{\star}\,:\,H^{1,\,1}_{BC}(X,\,\R)\to H^2_{DR}(X,\,\R), \hspace{2ex} [\alpha]_{BC}\mapsto \{\alpha\}_{DR},\end{equation}

\noindent is injective. (Of course, $P^{\star}$ is canonically well defined for any $X$ but it need not be injective if $X$ is not sGG.)

\end{Cor}

 That $P^{\star}$ is indeed the dual map of $P$ follows immediately from the identity $\int_X\alpha\wedge\Omega = \int_X\alpha\wedge\Omega^{n-1,\,n-1}$ which holds for bidegree reasons for any class $[\alpha]_{BC}\in H^{1,\,1}_{BC}(X,\,\R)$) and any class $\{\Omega\}_{DR}\in H^{2n-2}_{DR}(X,\,\R)$.

 No canonical right inverse of $P$ need exist when $X$ is only an sGG manifold, but for any given Hermitian metric on $X$ we will now construct a right inverse of $P$ depending on the chosen metric.

\begin{Def}\label{Def:Q-def} Let $X$ be an sGG compact complex manifold with $\mbox{dim}_{\C}X=n$ and let $\omega$ be an arbitrary Hermitian metric on $X$. With $\omega$ we associate the following injective linear map

\begin{equation}\label{eqn:Q-def}Q_{\omega}\,:\, H^{n-1,\,n-1}_A(X,\,\R)\to H^{2n-2}_{DR}(X,\,\R), \hspace{2ex} [\Omega^{n-1,\,n-1}]_A\mapsto \{\Omega\}_{DR},\end{equation}

\noindent where the real $d$-closed $(2n-2)$-form $\Omega$ on $X$ is determined by a given real $\partial\bar\partial$-closed $(n-1,\,n-1)$-form $\Omega^{n-1,\,n-1}$ and by the metric $\omega$ in the following way. \\

\noindent $(i)$\, If $\Delta_A$ denotes the Aeppli Laplacian associated with $\omega$, we have an orthogonal (w.r.t. the $L^2$ inner product defined by $\omega$) splitting

\vspace{1ex}

\hspace{10ex} $\ker(\partial\bar\partial) = \ker\Delta_A \oplus (\mbox{Im}\,\partial + \mbox{Im}\,\bar\partial)$  \hspace{6ex} (see e.g. [Pop13, $\S.2$]),

\vspace{1ex}

\noindent which induces a splitting of $\Omega^{n-1,\,n-1}\in\ker(\partial\bar\partial)$ as

\begin{equation}\label{eqn:Omega-splitting}\Omega^{n-1,\,n-1} = \Omega^{n-1,\,n-1}_A + \partial\Gamma^{n-2,\,n-1} + \bar\partial\Gamma^{n-1,\,n-2},\end{equation}

\noindent where $\Delta_A\Omega^{n-1,\,n-1}_A = 0$. The forms $\Gamma^{n-2,\,n-1}$, $\Gamma^{n-1,\,n-2}$ are of the shown bidegrees and are not uniquely determined, but we will see that the definition of $Q_{\omega}$ does not depend on their choices. (Since $\Omega^{n-1,\,n-1}$ is real, we can always choose $\Gamma^{n-2,\,n-1}=\overline{\Gamma^{n-1,\,n-2}}$.)

\vspace{1ex}

\noindent $(ii)$\, Since $\bar\partial(\partial\Omega^{n-1,\,n-1}_A) = 0$, we also have $d(\partial\Omega^{n-1,\,n-1}_A) = 0$. Thus the sGG assumption on $X$ and Lemma \ref{Lem:sGG-equiv-descriptions1} ensure that $\partial\Omega^{n-1,\,n-1}_A$ is $\bar\partial$-exact, i.e. there exists a smooth $(n,\,n-2)$-form $\Omega^{n,\,n-2}_A$ such that

\begin{equation}\label{eqn:d-bar-eqn-nn-2}\partial\Omega^{n-1,\,n-1}_A = \bar\partial(-\Omega^{n,\,n-2}_A).\end{equation}

\noindent We choose $\Omega^{n,\,n-2}_A$ to be the solution of equation (\ref{eqn:d-bar-eqn-nn-2}) of minimal $L^2$-norm (defined by $\omega$). Thus, $\Omega^{n,\,n-2}_A$ is uniquely determined by the formula

\begin{equation}\label{eqn:min-sol-nn-2}\Omega^{n,\,n-2}_A = -\bar\partial^{\star}\Delta^{''-1}(\partial\Omega^{n-1,\,n-1}_A),\end{equation}

\noindent where $\Delta''=\bar\partial\bar\partial^\star + \bar\partial^\star\bar\partial$ is the $\bar\partial$-Laplacian associated with the given metric $\omega$ and $\Delta^{''-1}$ is its Green operator (i.e. the inverse of its restriction to the orthogonal complement of its kernel).

\vspace{1ex}

\noindent $(iii)$\, Taking $\partial$ in (\ref{eqn:Omega-splitting}) and using (\ref{eqn:d-bar-eqn-nn-2}), we get\!\!:

$$\partial\Omega^{n-1,\,n-1} = \partial\Omega^{n-1,\,n-1}_A + \partial\bar\partial\Gamma^{n-1,\,n-2} = \bar\partial(-\Omega^{n,\,n-2}_A - \partial\Gamma^{n-1,\,n-2}).$$

\noindent We set \hspace{1ex} $\Omega^{n,\,n-2}:=\Omega^{n,\,n-2}_A + \partial\Gamma^{n-1,\,n-2}$. Thus we get\!\!:

\begin{equation}\label{eqn:partial-middle-dbar}\partial\Omega^{n-1,\,n-1} =  \bar\partial(-\Omega^{n,\,n-2}), \hspace{2ex}\mbox{hende also}\hspace{2ex} \bar\partial\Omega^{n-1,\,n-1} = \partial(-\overline{\Omega^{n,\,n-2}}),\end{equation}

\noindent where the latter identity follows from the former by taking conjugates and using the fact that $\Omega^{n-1,\,n-1}$ is real.

\vspace{1ex}

\noindent $(iv)$\, We set \hspace{1ex} $\Omega:= \Omega^{n,\,n-2} + \Omega^{n-1,\,n-1} + \overline{\Omega^{n,\,n-2}}$. It is clear that $\Omega$ is a real $(2n-2)$-form on $X$ and $d\Omega = 0$ (compare (\ref{eqn:partial-middle-dbar}) with (\ref{eqn:d-closedness-cond})).

\end{Def}

 We now make the trivial observation that Definition \ref{Def:Q-def} is correct.

\begin{Lem}\label{Lem:Q-def-correct} The map $Q_{\omega}$ is well defined and injective. Moreover, the composed linear map $P\circ Q_{\omega}\,:\,H^{n-1,\,n-1}_A(X,\,\R)\to H^{n-1,\,n-1}_A(X,\,\R)$ is the identity map of $H^{n-1,\,n-1}_A(X,\,\R)$ (so $Q_{\omega}$ is a section of $P$).

\end{Lem}

\noindent {\it Proof.} For well-definedness, we need to show that $Q_{\omega}([\Omega^{n-1,\,n-1}]_A)$ does not depend either on the choice of representative of the Aeppli class $[\Omega^{n-1,\,n-1}]_A$ or on the choice of the forms $\Gamma^{n-2,\,n-1}, \Gamma^{n-1,\,n-2}$ in (\ref{eqn:Omega-splitting}). Let us consider two {\it real} representatives of a same {\it real} Aeppli class\!\!:

$$[\Omega^{n-1,\,n-1}_1]_A = [\Omega^{n-1,\,n-1}_2]_A.$$

\noindent Let $\Omega_j=\Omega_j^{n,\,n-2} + \Omega_j^{n-1,\,n-1} + \overline{\Omega_j^{n,\,n-2}}$ ($j=1, 2$) be the real $d$-closed $(2n-2)$-forms on $X$ determined by $\Omega^{n-1,\,n-1}_j$ and $\omega$ as described in Definition \ref{Def:Q-def}.

 Since the $\Delta_A$-harmonic representative of a given Aeppli class is unique, we infer that $\Omega^{n-1,\,n-1}_{1,\,A} = \Omega^{n-1,\,n-1}_{2,\,A}$ (i.e. $\Omega^{n-1,\,n-1}_1$ and $\Omega^{n-1,\,n-1}_2$ have the same $\Delta_A$-harmonic projection). This implies that $\Omega^{n,\,n-2}_{1,\,A} = \Omega^{n,\,n-2}_{2,\,A}$ since the solution of minimal $L^2$-norm of a $\bar\partial$-equation (equation (\ref{eqn:d-bar-eqn-nn-2}) here) is unique. This further implies that

\begin{equation}\label{n_n-2_diff}\Omega^{n,\,n-2}_1 - \Omega^{n,\,n-2}_2 = \partial(\Gamma^{n-1,\,n-2}_1 - \Gamma^{n-1,\,n-2}_2).\end{equation}

\noindent On the other hand, (\ref{eqn:Omega-splitting}) spells

$$\Omega^{n-1,\,n-1}_j = \Omega^{n-1,\,n-1}_{j,\,A} + \partial\Gamma^{n-2,\,n-1}_j + \bar\partial\Gamma^{n-1,\,n-2}_j, \hspace{3ex} j=1, 2,$$

\noindent which gives, since $\Omega^{n-1,\,n-1}_{1,\,A} = \Omega^{n-1,\,n-1}_{2,\,A}$, the identity

\begin{equation}\label{n-1_n-1_diff}\Omega^{n-1,\,n-1}_1 - \Omega^{n-1,\,n-1}_2 = \partial(\Gamma^{n-2,\,n-1}_1 - \Gamma^{n-2,\,n-1}_2) + \bar\partial(\Gamma^{n-1,\,n-2}_1 - \Gamma^{n-1,\,n-2}_2).\end{equation}

\noindent We see that (\ref{n_n-2_diff}) and (\ref{n-1_n-1_diff}) amount to the $d$-exactness condition (\ref{eqn:d-exactness-cond}) for the real $(2n-2)$-form $\Omega:=\Omega_1 - \Omega_2$ (where we choose $\Gamma^{n,\,n-3}=0$ and $\Gamma^{n-3,\,n}=0$). Thus $\Omega_1 - \Omega_2$ is $d$-exact, i.e. $\{\Omega_1\}_{DR} = \{\Omega_2\}_{DR}$, so $Q_{\omega}([\Omega^{n-1,\,n-1}_1]_A) = Q_{\omega}([\Omega^{n-1,\,n-1}_2]_A)$.

 To show that $Q_{\omega}$ is injective, let $Q_{\omega}([\Omega^{n-1,\,n-1}]_A)=0$, i.e. $\{\Omega\}_{DR} = 0$. This means that $\Omega$ is $d$-exact, which in turn means that the identities (\ref{eqn:d-exactness-cond}) hold. It is clear that $(ii)$ of (\ref{eqn:d-exactness-cond}) expresses the fact that $[\Omega^{n-1,\,n-1}]_A=0$.

 The fact that $P\circ Q_{\omega}([\Omega^{n-1,\,n-1}]_A) = [\Omega^{n-1,\,n-1}]_A$ for any class $[\Omega^{n-1,\,n-1}]_A\in H^{n-1,\,n-1}_A(X,\,\R)$ follows immediately from the definitions of $P$ and $Q_{\omega}$\!\!: the original $(n-1,\,n-1)$-form $\Omega^{n-1,\,n-1}$ is indeed the $(n-1,\,n-1)$-component of the $(2n-2)$-form constructed from $\Omega^{n-1,\,n-1}$ in Definition \ref{Def:Q-def}.  \hfill $\Box$

\vspace{2ex}

 Putting these pieces of information together, we immediately get the

\begin{Cor}\label{Cor:Qstar-surj} Let $X$ be an sGG compact complex manifold, $\mbox{dim}_{\C}X=n$. For any Hermitian metric $\omega$ on $X$, the dual map of $Q_{\omega}$\!\!:

\begin{equation}\label{eqn:Pstar-inj} Q_{\omega}^{\star}\,:\,H^2_{DR}(X,\,\R)\to H^{1,\,1}_{BC}(X,\,\R)\end{equation}

\noindent is surjective. Moreover, the composition $Q_{\omega}^{\star}\circ P^{\star}\,:\,H^{1,\,1}_{BC}(X,\,\R)\to H^{1,\,1}_{BC}(X,\,\R)$ is the identity map.

\end{Cor}

 Note that the dual map $Q_{\omega}^{\star}$ has the following explicit form

\begin{equation}\label{eqn:Qstar-def}\int\limits_XQ_{\omega}^{\star}(\{\alpha\}_{DR})\wedge[\Omega^{n-1,\,n-1}]_A = \int\limits_X\{\alpha\}_{DR}\wedge Q_{\omega}([\Omega^{n-1,\,n-1}]_A)\end{equation}

\noindent for any classes $\{\alpha\}_{DR}\in H^2_{DR}(X,\,\R)$ and $[\Omega^{n-1,\,n-1}]_A\in H^{n-1,\,n-1}_A(X,\,\R)$. (The meaning of cohomology classes in (\ref{eqn:Qstar-def}) is that the integrals do not depend on the choice of representatives in those classes.)

\vspace{2ex}

 We can well call the pair of maps $(P,\, Q_{\omega})$ a {\it fake Hodge-Aeppli decomposition} of $H^{2n-2}_{DR}(X,\,\R)$ and the dual pair of maps $(P^{\star},\,Q_{\omega}^{\star})$ the {\it dual fake Hodge-Bott-Chern decomposition} of $H^2_{DR}(X,\,\R)$.

\subsection{Deformation semicontinuity of ${\cal G}_X$ and ${\cal E}_X$}\label{subsection:s-c}

 We now use the fake Hodge decomposition of the previous subsection to study the variations of the Gauduchon cone ${\cal G}_X$ and of the pseudo-effective cone ${\cal E}_X$ under small deformations of the sGG complex structure of $X$.

 Let $\pi\,:\,{\cal X}\to\Delta$ be a holomorphic family of compact complex manifolds. Without loss of generality, we may suppose that $\Delta\subset\C$ is an open disc about the origin. The fibres $X_t:=\pi^{-1}(t)\subset{\cal X}$ ($t\in\Delta$) are thus compact complex manifolds of equal dimension $n$ and the family is $C^{\infty}$ locally trivial, hence the De Rham cohomology groups of the fibres can be identified with a fixed space $H^k(X,\,\C)$ for every $k=0, 1, \dots , 2n$. As the complex structure of $X_t$ varies with $t\in\Delta$, the Bott-Chern, Dolbeault and Aeppli cohomologies of the fibres depend on $t$.

 Suppose moreover that $X_0$ is an sGG manifold. Then $X_t$ is an sGG manifold for all $t\in\Delta$ sufficiently close to $0$ by our Corollary \ref{Cor:openness}. After shrinking $\Delta$ about $0$, we can assume that $X_t$ is an sGG manifold for all $t\in\Delta$. We fix any $C^{\infty}$ family $(\omega_t)_{t\in\Delta}$ of Hermitian metrics on the respective fibres $(X_t)_{t\in\Delta}$. Let $t_0\in\Delta$ be an arbitrary point (e.g. we take $t_0=0$).

\vspace{2ex}

\noindent $(I)$\, {\it Variation of the Gauduchon cone} \\

 The {\it fake Hodge-Aeppli decomposition} constructed in the previous subsection on each fibre $X_t$ gives us maps as follows\!\!:

\vspace{2ex}

$\begin{array}{llllll}\label{fHA-G} H_A^{n-1,\,n-1}(X_0,\,\R) & \stackrel{Q_{\omega_0}}{\hookrightarrow} & H^{2n-2}_{DR}(X,\,\R) & \stackrel{P_t}{\twoheadrightarrow} & H_A^{n-1,\,n-1}(X_t,\,\R), &  t\in\Delta. \\
    \bigcup &   &   &   &  \bigcup\\
    {\cal G}_{X_0} & & & &   {\cal G}_{X_t} &   \end{array}$

\vspace{2ex}

\noindent Thus the image of the Gauduchon cone ${\cal G}_{X_0}$ of $X_0$ under the composition $P_t\circ Q_{\omega_0}\,:\,H_A^{n-1,\,n-1}(X_0,\,\R)\to H_A^{n-1,\,n-1}(X_t,\,\R)$ can be compared to ${\cal G}_{X_t}$ as subsets of $H_A^{n-1,\,n-1}(X_t,\,\R)$. Note that $P_0\circ Q_{\omega_0} = \mbox{Id}_{H_A^{n-1,\,n-1}(X_0,\,\R)}$ and it follows  from [KS60, Theorem 5] that the surjections $(P_t)_{t\in\Delta}$ vary in a $C^{\infty}$ way with $t$ (hence the maps $P_t\circ Q_{\omega_0}$ are isomorphisms) if the dimension of $H_A^{n-1,\,n-1}(X_t,\,\R)$ ($=h_{BC}^{1,\,1}(t)$ by duality) is independent of $t$. However, if $h_{BC}^{1,\,1}(0) > h_{BC}^{1,\,1}(t)$ for $t\neq 0$ close to $0$, the Gauduchon cone ${\cal G}_{X_0}$ of $X_0$ has more dimensions than its counterparts ${\cal G}_{X_t}$ on the nearby fibres, but the projections $P_t$ for $t\neq 0$ eliminate the extra dimensions of $(P_t\circ Q_{\omega_0})({\cal G}_{X_0})$. It seems sensible to introduce the following definition.

\begin{Def}\label{Def:lim-G} If $(X_t)_{t\in\Delta}$ is a holomorphic family of {\bf sGG} compact complex $n$-dimensional manifolds, the {\bf limit as $t\rightarrow t_0$ of the Gauduchon cones ${\cal G}_{X_t}$} of the fibres $X_t$ for $t\neq t_0$ is defined as the following subset of $H_A^{n-1,\,n-1}(X_{t_0},\,\R)$\!\!:

$$\lim\limits_{t\rightarrow t_0}{\cal G}_{X_t}:=\bigg\{[\Omega^{n-1,\,n-1}]_A\in H_A^{n-1,\,n-1}(X_{t_0},\,\R)\,\mid\,(P_t\circ Q_{\omega_{t_0}})([\Omega^{n-1,\,n-1}]_A)\in{\cal G}_{X_t}\hspace{1ex}\forall t\sim t_0\bigg\},$$

\noindent where ``\,$\forall t\sim t_0$'' means ``for all $t$ sufficiently close to $t_0$''.

\end{Def}

 Note that $\lim\limits_{t\rightarrow t_0}{\cal G}_{X_t}$ depends on the metric $\omega_{t_0}$ (since $Q_{\omega_{t_0}}$ depends thereon) but does not depend on the way in which $\omega_{t_0}$ has been extended in a $C^{\infty}$ fashion to metrics $\omega_t$ on the nearby fibres (since the maps $P_t$ are canonical).

 We can now prove that the Gauduchon cone ${\cal G}_{X_t}$ of the sGG fibre $X_t$ behaves lower-semicontinuously w.r.t. $t\in\Delta$ in the usual topology of $\Delta$ much as it did in the special case of families of $\partial\bar\partial$-manifolds treated in [Pop13b].

\begin{The}\label{The:lsc-G} Let $(X_t)_{t\in\Delta}$ be any holomorphic family of {\bf sGG} compact complex manifolds endowed with any $C^{\infty}$ family $(\omega_t)_{t\in\Delta}$ of Hermitian metrics.

 Then, for all $t_0\in\Delta$, the following inclusion holds\!\!:

$${\cal G}_{X_{t_0}}\subset\lim\limits_{t\rightarrow t_0}{\cal G}_{X_t}.$$

\end{The}

\noindent {\it Proof.} We may assume that $t_0=0$. Denote by $n$ the complex dimension of the fibres and let $[\gamma_0^{n-1}]_A\in{\cal G}_{X_0}$ for some Gauduchon metric $\gamma_0>0$ on $X_0$. Let $\Omega$ be the $C^{\infty}$ real $d$-closed $(2n-2)$-form determined by $\gamma_0^{n-1}$ and by the Hermitian metric $\omega_0$ as in Definition \ref{Def:Q-def}. For every $t\in\Delta$, the splitting of $\Omega$ into pure-type forms reads\!\!:

$$\Omega = \Omega_t^{n,\,n-2} + \Omega_t^{n-1,\,n-1} + \Omega_t^{n-2,\,n}, \hspace{2ex} t\in\Delta,$$

\noindent where $\Omega_0^{n-1,\,n-1} = \gamma_0^{n-1}$. Then $(P_t\circ Q_{\omega_0})([\gamma_0^{n-1}]_A) = P_t(\{\Omega\}_{DR}) = [\Omega_t^{n-1,\,n-1}]_A\in H^{n-1,\,n-1}_A(X_t,\,\R)$ for all $t\in\Delta$. Since $\Omega_0^{n-1,\,n-1}>0$ and the $\Omega_t^{n-1,\,n-1}$ vary in a $C^{\infty}$ way with $t\in\Delta$ (as components of the fixed form $\Omega$), we get

$$\Omega_t^{n-1,\,n-1}>0 \hspace{2ex}\mbox{for all $t$ sufficiently close to $0$},$$

\noindent hence there exists a unique Gauduchon metric $\gamma_t$ on $X_t$ such that $\gamma_t^{n-1}=\Omega_t^{n-1,\,n-1}$, so $[\Omega_t^{n-1,\,n-1}]_A\in{\cal G}_{X_t}$ for all $t$ close to $0$. Thus $[\gamma_0^{n-1}]_A\in\lim\limits_{t\rightarrow t_0}{\cal G}_{X_t}$.

\hfill $\Box$

\vspace{2ex}

\noindent $(II)$\, {\it Dual situation\!: variation of the pseudo-effective cone} \\

 The dual of the {\it fake Hodge-Aeppli decomposition} constructed in the previous subsection on each fibre $X_t$ gives us maps as follows\!\!:

\vspace{2ex}

$\begin{array}{llllll} H_{BC}^{1,\,1}(X_0,\,\R) & \stackrel{P_0^{\star}}{\hookrightarrow} & H^2_{DR}(X,\,\R) & \stackrel{Q_{\omega_t}^{\star}}{\twoheadrightarrow} & H_{BC}^{1,\,1}(X_t,\,\R), & t\in\Delta.\\
    \bigcup &   &   &   &  \bigcup\\
    {\cal E}_{X_0} & & & &   {\cal E}_{X_t} &  \end{array}$

\vspace{2ex}

\noindent Clearly, $Q_{\omega_0}^{\star}\circ P_0^{\star} = \mbox{Id}_{H_{BC}^{1,\,1}(X_0,\,\R)}$ and it follows from [KS60, Theorem 5] that the surjections $(Q_{\omega_t}^{\star})_{t\in\Delta}$ vary in a $C^{\infty}$ way with $t$ (hence the maps $Q_{\omega_t}^{\star}\circ P_0^{\star}$ are isomorphisms) if the dimension of $H_{BC}^{1,\,1}(X_t,\,\R)$ ($=h_{BC}^{1,\,1}(t)$) is independent of $t$. However, if $h_{BC}^{1,\,1}(0) > h_{BC}^{1,\,1}(t)$ for $t\neq 0$ close to $0$, the pseudo-effective cone ${\cal E}_{X_0}$ of $X_0$ lies in a space which has a higher dimension than the ambient spaces of its counterparts ${\cal E}_{X_t}$ on the nearby fibres.

 We shall now define a complex vector subspace

$$H_{BC}^{'1,\,1}(X_0,\,\C)\subset H_{BC}^{1,\,1}(X_0,\,\C)$$

\noindent depending on the chosen family of Hermitian metrics $(\omega_t)_{t\in\Delta}$ such that\!\!:

 $\cdot$\, $\mbox{dim}\,H_{BC}^{'1,\,1}(X_0,\,\C) = \mbox{dim}\,H_{BC}^{1,\,1}(X_t,\,\C)$ for all $t\in\Delta$ close to $0$\!;

 $\cdot$\, $H_{BC}^{'1,\,1}(X_0,\,\C) = H_{BC}^{1,\,1}(X_0,\,\C)$ when $h_{BC}^{1,\,1}(0) = h_{BC}^{1,\,1}(t)$ for $t$ close to $0$.

\vspace{2ex}

\noindent For every $t\in\Delta$, let $\Delta_{A,\,t}\,:\,C^{\infty}_{n-1,\,n-1}(X_t,\,\C)\to C^{\infty}_{n-1,\,n-1}(X_t,\,\C)$ be the Aeppli Laplacian in bidegree $(n-1,\,n-1)$ defined by the Hermitian metric $\omega_t$ on $X_t$. Since $\Delta_{A,\,t}$ is a non-negative self-adjoint elliptic operator (of order $4$), it has a discrete spectrum $0\leq \lambda_1(t)\leq \lambda_2(t)\leq \dots  $ with $+\infty$ as sole accumulation point and the space $C^{\infty}_{n-1,\,n-1}(X_t,\,\C)$ has an orthonormal basis $(e_j(t))_{j\geq 1}$ consisting of eigenvectors such that

$$\Delta_{A,\,t}e_j(t) = \lambda_j(t)\,e_j(t), \hspace{2ex} j\geq 1, t\in\Delta.$$

\noindent Let $N:=h_{BC}^{1,\,1}(0)$ and $p:=h_{BC}^{1,\,1}(t)$ for $t$ close to $0$. Thus $N\geq p$ by the Kodaira-Spencer upper-semicontinuity property [KS60, Theorem 4]. Let

\vspace{1ex}

\noindent$0 < \varepsilon < \min\bigg(\mbox{Spec}\,\Delta_{A,\,0}\cap(0,\,+\infty)\bigg) \hspace{2ex} \mbox{such that}\hspace{2ex}  \varepsilon\notin\mbox{Spec}\,\Delta_{A,\,t} \hspace{1ex} \forall t\sim 0.$

\vspace{1ex}

\noindent Then, thanks to fundamental Kodaira-Spencer theorems on smooth families of elliptic operators [KS60, Theorems 1-5], we have the following picture\!\!: \\

\noindent $0 = \lambda_1(0) = \dots = \lambda_N(0) < \varepsilon < \lambda_{N+1}(0)$, while for all $t\sim 0, t\neq 0$, we have\!\!:

\noindent $0 = \lambda_1(t) = \dots = \lambda_p(t) < \lambda_{p+1}(t) \leq \dots \leq \lambda_N(t) < \varepsilon < \lambda_{N+1}(t),$ \\

\noindent i.e. the number of eigenvalues (counted with multiplicities) of $\Delta_{A,\,t}$ lying in the open interval $(-1,\,\varepsilon)$ is independent of $t$ if $t\in\Delta$ is sufficiently close to $0$. Moreover, if $E_{\Delta_{A,\,t}}(\lambda)$ denotes the eigenspace of $\Delta_{A,\,t}$ corresponding to the eigenvalue $\lambda$, the Kodaira-Spencer theorems further ensure that

$$\Delta\ni t\mapsto \bigoplus\limits_{0\leq\lambda<\varepsilon}E_{\Delta_{A,\,t}}(\lambda):={\cal E}_{A,\,\varepsilon}^{n-1,\,n-1}(t)$$

\noindent is a $C^{\infty}$ vector bundle of finite rank (equal to $N$ here) after possibly shrinking $\Delta$ about $0$ and that the orthogonal projections

\begin{equation}\label{eqn:sigma_t-def}C^{\infty}_{n-1,\,n-1}(X_t,\,\C)\stackrel{\sigma_t}{\longrightarrow}{\cal E}_{A,\,\varepsilon}^{n-1,\,n-1}(t)\end{equation}

\noindent vary in a $C^{\infty}$ way with $t\in\Delta$.

\noindent Thus $\{e_1(t), \dots , e_N(t)\}$ is a local frame for the vector bundle ${\cal E}_{A,\,\varepsilon}^{n-1,\,n-1}$ and we have\!\!:

\begin{eqnarray}\nonumber\ker\Delta_{A,\,0} & = & \langle e_1(0), \dots , e_p(0),\dots, e_N(0)\rangle = {\cal E}_{A,\,\varepsilon}^{n-1,\,n-1}(0),\\
\nonumber  \ker\Delta_{A,\,t} & = & \langle e_1(t), \dots , e_p(t)\rangle\subset\langle e_1(t), \dots , e_N(t)\rangle = {\cal E}_{A,\,\varepsilon}^{n-1,\,n-1}(t), \hspace{2ex} t\neq 0.\end{eqnarray}

\noindent Thus we have an orthogonal splitting

\vspace{1ex}

$\ker\Delta_{A,\,0} = \langle e_1(0), \dots , e_p(0)\rangle \oplus \langle e_{p+1}(0),\dots , e_N(0)\rangle$

\vspace{1ex}

\noindent which induces under the Hodge isomorphism $\ker\Delta_{A,\,0}\simeq H^{n-1,\,n-1}_A(X_0,\,\C)$ a splitting

\begin{equation}\label{eqn:H_a-splitting}H^{n-1,\,n-1}_A(X_0,\,\C) = H^{'n-1,\,n-1}_A(X_0,\,\C)\oplus H^{''n-1,\,n-1}_A(X_0,\,\C)\end{equation}

\noindent where $H^{'n-1,\,n-1}_A(X_0,\,\C)\simeq\langle e_1(0), \dots , e_p(0)\rangle$ and

\noindent $H^{''n-1,\,n-1}_A(X_0,\,\C)\simeq\langle e_{p+1}(0), \dots , e_N(0)\rangle$. Now, $H^{1,\,1}_{BC}(X_0,\,\C)$ and $H^{n-1,\,n-1}_A(X_0,\,\C)$ are dual to each other, so identifying $H^{1,\,1}_{BC}(X_0,\,\C)$ with $H^{n-1,\,n-1}_A(X_0,\,\C)^{\star}$ we define

\begin{equation}\label{eqn:H11_BC'-def}H^{'1,\,1}_{BC}(X_0,\,\C):=\bigg\{[\alpha]_{BC}\in H^{1,\,1}_{BC}(X_0,\,\C)\,\mid\, [\alpha]_{BC|H^{''n-1,\,n-1}_A(X_0,\,\C)} = 0\bigg\}.\end{equation}

\noindent Thus, $H^{'1,\,1}_{BC}(X_0,\,\C)$ consists of the linear maps $[\alpha]_{BC}\,:\,H^{n-1,\,n-1}_A(X_0,\,\C)\to\C$ vanishing on $H^{''n-1,\,n-1}_A(X_0,\,\C)$, i.e. identifies with the dual of $H^{'n-1,\,n-1}_A(X_0,\,\C)$.

 It is clear that $H^{'1,\,1}_{BC}(X_0,\,\C)$ coincides with $H^{1,\,1}_{BC}(X_0,\,\C)$ if $h^{1,\,1}_{BC}(0) = h^{1,\,1}_{BC}(t)$ for $t\sim 0$, but it depends on the choice of the $C^{\infty}$ family of metrics $(\omega_t)_{t\in\Delta}$, so it is not canonical, if $h^{1,\,1}_{BC}(0) > h^{1,\,1}_{BC}(t)$. The same construction can, of course, be run for any $t_0\in\Delta$ in place of $0$.

\begin{Def}\label{Def:lim-E} Let $(X_t)_{t\in\Delta}$ be a holomorphic family of {\bf sGG} compact complex manifolds equipped with a $C^{\infty}$ family of Hermitian metrics $(\omega_t)_{t\in\Delta}$.

 For any $t_0\in\Delta$, the {\bf limit as $t\rightarrow t_0$ of the pseudo-effective cones ${\cal E}_{X_t}$} of the fibres $X_t$ for $t\neq t_0$ is defined as the following subset of $H_{BC}^{1,\,1}(X_{t_0},\,\R)$\!\!:

$$\lim\limits_{t\rightarrow t_0}{\cal E}_{X_t}:=\bigg\{[\alpha]_{BC}\in H_{BC}^{1,\,1}(X_{t_0},\,\R)\cap H_{BC}^{'1,\,1}(X_{t_0},\,\C)\,\mid\,(Q_{\omega_t}^{\star}\circ P_{t_0}^{\star})([\alpha]_{BC})\in{\cal E}_{X_t}\hspace{1ex}\forall t\sim t_0\bigg\},$$

\noindent where ``\,$\forall t\sim t_0$'' means ``for all $t$ sufficiently close to $t_0$''.

\end{Def}

 Note that we restrict from the start to classes in the subspace $H_{BC}^{'1,\,1}(X_{t_0},\,\C)$ $\subset H_{BC}^{1,\,1}(X_{t_0},\,\C)$ to trim off the extra dimensions that the limit may acquire if the dimension of $H_{BC}^{1,\,1}(X_t,\,\C)$ increases in the limit. Note also that, unlike its Gauduchon-cone counterpart, $\lim\limits_{t\rightarrow t_0}{\cal E}_{X_t}$ depends not only on the metric $\omega_{t_0}$ but on the whole family of metrics $(\omega_t)_{t\in\Delta}$ for $t\sim 0$.

 We can now prove that the pseudo-effective cone ${\cal E}_{X_t}$ behaves upper-semicontinuously in families of sGG manifolds.

\begin{The}\label{The:usc-E} Let $(X_t)_{t\in\Delta}$ be any holomorphic family of {\bf sGG} compact complex manifolds endowed with any $C^{\infty}$ family $(\omega_t)_{t\in\Delta}$ of Hermitian metrics.

 Then, for all $t_0\in\Delta$, the following inclusion holds\!\!:

$${\cal E}_{X_{t_0}}\supset\lim\limits_{t\rightarrow t_0}{\cal E}_{X_t}.$$

\end{The}

\noindent {\it Proof.} Without loss of generality, we may suppose that $t_0=0$. Let $[T]_{BC}\in\lim\limits_{t\rightarrow 0}{\cal E}_{X_t}$, where $T$ is a $d$-closed real $(1,\,1)$-current on $X_0$. (We implicitly use the fact that the Bott-Chern cohomology can be computed using either smooth forms or currents.) Since ${\cal E}_{X_t}$ is the dual of $\overline{{\cal G}_{X_t}}$ by Lamari's duality lemma, for all $t\neq 0$ with $t\sim 0$ and for any Gauduchon metric $\gamma_t$ on $X_t$, we have

\begin{eqnarray}\label{eqn:integral-choiceT}\int\limits_{X_t}(Q_{\omega_t}^{\star}\circ P_0^{\star})([T]_{BC})\wedge[\gamma_t^{n-1}]_A & = & \int\limits_{X_t}\{T\}_{DR}\wedge Q_{\omega_t}([\gamma_t^{n-1}]_A)\\
\nonumber & = & \int\limits_{X_t}T\wedge(\Omega_t^{n,\,n-2} + \gamma_t^{n-1} + \overline{\Omega_t^{n,\,n-2}})\geq 0.\end{eqnarray}

\noindent Indeed, the first identity in (\ref{eqn:integral-choiceT}) is (\ref{eqn:Qstar-def}), while the second identity holds for the $(n,\,n-2)$-form $\Omega_t^{n,\,n-2}$ on $X_t$ determined as described in Definition \ref{Def:Q-def} by $\gamma_t^{n-1}$ and the Hermitian metric $\omega_t$ of $X_t$.

 We will show that $[T]_{BC}\in{\cal E}_{X_0}$. Since ${\cal E}_{X_0}$ is the dual of $\overline{{\cal G}_{X_0}}$ by Lamari's lemma, this amounts to showing that

\begin{equation}\label{eqn:T_to-show}\int\limits_{X_0}T\wedge\gamma_0^{n-1}\geq 0\end{equation}

\noindent for any Gauduchon metric $\gamma_0$ on $X_0$.

 Let us fix an arbitrary Gauduchon metric $\gamma_0$ on $X_0$. Pick any $C^{\infty}$ deformation of $\gamma_0$ to Gauduchon metrics $(\gamma_t)_{t\in\Delta}$ on the fibres $(X_t)_{t\in\Delta}$. (This is always possible as the proof of Gauduchon's theorem shows -- see e.g. [Pop13a, $\S.3$]). Since $(\gamma_t^{n-1})_{t\in\Delta}$ is a $C^{\infty}$ family of $(n-1,\,n-1)$-forms and since $(\sigma_t)_{t\in\Delta}$ is a $C^{\infty}$ family of orthogonal projections (defined in (\ref{eqn:sigma_t-def})), $(\sigma_t\gamma_t^{n-1})_{t\in\Delta}$ is a $C^{\infty}$ family of $(n-1,\,n-1)$-forms.

 We use the notation of Definition \ref{Def:Q-def} with $\Omega^{n-1,\,n-1}$ replaced with $\gamma_t^{n-1}$ on each fibre $X_t$. Thus $\Omega^{n-1,\,n-1}_{A,\,t}$ stands for the $\Delta_{A,\,t}$-harmonic component of $\gamma_t^{n-1}$, so for every $t\in\Delta$ we have\!\!:

\vspace{1ex}

$\gamma_t^{n-1} = \Omega^{n-1,\,n-1}_{A,\,t} + \partial_t\Gamma_t^{n-2,\,n-1} + \bar\partial_t\overline{\Gamma_t^{n-2,\,n-1}},$

\vspace{1ex}

$\Omega^{n,\,n-2}_{A,\,t}:= - \bar\partial_t^{\star}\Delta_t^{''-1}(\partial_t\Omega^{n-1,\,n-1}_{A,\,t})$ \hspace{2ex} and \hspace{2ex} $\Omega^{n,\,n-2}_t:= \Omega^{n,\,n-2}_{A,\,t} + \partial_t\overline{\Gamma_t^{n-2,\,n-1}}$

\vspace{1ex}

\noindent On the other hand, we have\!\!:

\vspace{1ex}

$\sigma_0\gamma_0^{n-1} = \Omega^{n-1,\,n-1}_{A,\,0} = \sum\limits_{j=1}^pc_j(0)\,e_j(0) + \sum\limits_{j=p+1}^Nc_j(0)\,e_j(0) := \Omega'_{A,\,0} + \Omega''_{A,\,0},$

\vspace{1ex}

$\sigma_t\gamma_t^{n-1} = \sum\limits_{j=1}^pc_j(t)\,e_j(t) + \sum\limits_{j=p+1}^Nc_j(t)\,e_j(t) = \Omega^{n-1,\,n-1}_{A,\,t} + \sum\limits_{j=p+1}^Nc_j(t)\,e_j(t),$

\vspace{1ex}

\noindent for $t\sim 0, t\neq 0$. Thus $\Omega'_{A,\,0}, \Omega''_{A,\,0}\in\ker\Delta_{A,\,0}$ and $[\Omega'_{A,\,0}]_A\in H^{'n-1,\,n-1}_A(X_0,\,\C)$ while $[\Omega''_{A,\,0}]_A\in H^{''n-1,\,n-1}_A(X_0,\,\C)$. The coefficients $c_j(t)\in\C$ vary continuously with $t\in\Delta$, so $c_j(t)\rightarrow c_j(0)$ as $t\rightarrow 0$ for every $j$. We get\!\!:

\begin{equation}\label{eqn:Omega-n-1-conv}\Omega^{n-1,\,n-1}_{A,\,t} = \sum\limits_{j=1}^pc_j(t)\,e_j(t) \longrightarrow \sum\limits_{j=1}^pc_j(0)\,e_j(0) =  \Omega'_{A,\,0}\in\ker\Delta_{A,\,0}   \hspace{2ex} \mbox{as}\hspace{1ex} t\rightarrow 0,\end{equation}

\noindent hence, from $\partial_t$ varying in a $C^{\infty}$ way with $t$ up to $t=0$, we infer

$$\partial_t\Omega^{n-1,\,n-1}_{A,\,t} \longrightarrow \partial_0\Omega'_{A,\,0}  \hspace{2ex} \mbox{as}\hspace{1ex} t\rightarrow 0.$$

\noindent Now comes a crucial argument. The forms $\partial_t\Omega^{n-1,\,n-1}_{A,\,t}$ and $\partial_0\Omega'_{A,\,0}$ are of bidegree $(n,\,n-1)$ for their respective complex structures. On the other hand, by Serre duality we have $h^{n,\,n-1}_{\bar\partial}(t) = h^{0,\,1}_{\bar\partial}(t)$, hence part $(ii)$ of our Corollary \ref{Cor:openness} and the sGG assumption ensure that

$$h^{n,\,n-1}_{\bar\partial}(0) = h^{n,\,n-1}_{\bar\partial}(t) \hspace{2ex} \mbox{for all}\hspace{1ex}t\sim 0.$$

\noindent Therefore, the Green operators $(\Delta_t^{''-1})_{t\in\Delta}$ vary in a $C^{\infty}$ way with $t$ (up to $t=0$) by the Kodaira-Spencer theorem [KS60, Theorem 5] which applies when the relevant Hodge numbers ($h^{n,\,n-1}_{\bar\partial}(t)$ here) do not jump. Thus,

$$\Delta_t^{''-1}(\partial_t\Omega^{n-1,\,n-1}_{A,\,t}) \longrightarrow \Delta_0^{''-1}(\partial_0\Omega'_{A,\,0})  \hspace{2ex} \mbox{as}\hspace{1ex} t\rightarrow 0$$

\noindent and since $\bar\partial_t^{\star}$ varies in a $C^{\infty}$ way with $t$ up to $t=0$, we infer

\begin{equation}\label{eqn:nn-2_t-convergence}\Omega^{n,\,n-2}_{A,\,t}= - \bar\partial_t^{\star}\Delta_t^{''-1}(\partial_t\Omega^{n-1,\,n-1}_{A,\,t}) \longrightarrow  - \bar\partial_0^{\star}\Delta_0^{''-1}(\partial_0\Omega'_{A,\,0}):= \Omega^{'n,\,n-2}_{A,\,0}\end{equation}

\noindent as $t\rightarrow 0$. It is clear that the form $\Omega^{'n,\,n-2}_{A,\,0}$ is of bidegree $(n,\,n-2)$ on $X_0$.

 We can now finish the proof of the theorem. Recall that we have to prove inequality (\ref{eqn:T_to-show}). With the above preparations, we have\!\!:

\begin{eqnarray}\label{eqn:T_to-show-bis}\nonumber\int\limits_{X_0}T\wedge\gamma_0^{n-1} & = & \int\limits_{X_0}T\wedge(\Omega^{n-1,\,n-1}_{A,\,0} + \partial_0\Gamma_0^{n-2,\,n-1} + \bar\partial_0\overline{\Gamma_0^{n-2,\,n-1}}) \stackrel{(a)}{=} \int\limits_{X_0}T\wedge\Omega^{n-1,\,n-1}_{A,\,0}\\
    & = & \int\limits_{X_0}T\wedge\Omega'_{A,\,0} + \int\limits_{X_0}T\wedge\Omega''_{A,\,0} \stackrel{(b)}{=} \int\limits_{X_0}T\wedge\Omega'_{A,\,0},\end{eqnarray}

\noindent where identity $(a)$ follows by Stokes' theorem from $\partial_0T=0$ and $\bar\partial_0T=0$ (due to $dT=0$), while identity $(b)$ follows from the definition of $H^{'1,\,1}_{BC}(X_0,\,\C)$, from $[T]_{BC}\in H^{'1,\,1}_{BC}(X_0,\,\C)$ and from $[\Omega''_{A,\,0}]_A\in H^{''n-1,\,n-1}_A(X_0,\,\C)$.

 On the other hand, the last integral in (\ref{eqn:integral-choiceT}), which is non-negative for all $t\sim 0$ and $t\neq 0$, transforms as follows\!\!:

\begin{eqnarray}\label{eqn:1main-integral_t}\int\limits_{X_t}T\wedge(\Omega_t^{n,\,n-2} + \gamma_t^{n-1} + \overline{\Omega_t^{n,\,n-2}}) & = & \int\limits_{X_t}T\wedge\gamma_t^{n-1}\\
\nonumber & + & \int\limits_{X_t}T\wedge\Omega_{A,\,t}^{n,\,n-2} + \int\limits_{X_t}T\wedge\partial_t\Gamma_t^{n-1,\,n-2}\\
\nonumber   & + & \int\limits_{X_t}T\wedge\overline{\Omega_{A,\,t}^{n,\,n-2}} + \int\limits_{X_t}T\wedge\bar\partial_t\Gamma_t^{n-2,\,n-1}.\end{eqnarray}

\noindent Now, $\int_{X_t}T\wedge\gamma_t^{n-1}$ converges to $\int_{X_0}T\wedge\gamma_0^{n-1}$, while $\int_{X_t}T\wedge\Omega_{A,\,t}^{n,\,n-2}$ converges to $\int_{X_0}T\wedge\Omega_{A,\,0}^{'n,\,n-2} = 0$ by the crucial convergence (\ref{eqn:nn-2_t-convergence}). The last identity follows from $T$ being of bidegree $(1,\,1)$ and $\Omega_{A,\,0}^{'n,\,n-2}$ being of bidegree $(n,\,n-2)$, hence $T\wedge\Omega_{A,\,0}^{'n,\,n-2} = 0$ as an $(n+1,\, n-1)$-current. By conjugation, we infer that $\int_{X_t}T\wedge\overline{\Omega_{A,\,t}^{n,\,n-2}}$ converges to $\int_{X_0}T\wedge\overline{\Omega_{A,\,0}^{'n,\,n-2}} = 0$. Furthermore, we have\!\!:

$$\int\limits_{X_t}T\wedge\partial_t\overline{\Gamma_t^{n-2,\,n-1}} = \int\limits_{X_t}T\wedge d\overline{\Gamma_t^{n-2,\,n-1}} - \int\limits_{X_t}T\wedge\bar\partial_t\overline{\Gamma_t^{n-2,\,n-1}} =  - \int\limits_{X_t}T\wedge\bar\partial_t\overline{\Gamma_t^{n-2,\,n-1}},$$

\noindent the last identity following from Stokes' theorem and $dT=0$. We also have the conjugate identity\!\!: $\int_{X_t}T\wedge\bar\partial_t\Gamma_t^{n-2,\,n-1} = - \int_{X_t}T\wedge\partial_t\Gamma_t^{n-2,\,n-1}$, hence\!\!:

\begin{eqnarray}\label{eqn:2main-integral_t}\nonumber\int\limits_{X_t}T\wedge\partial_t\overline{\Gamma_t^{n-2,\,n-1}} + \int\limits_{X_t}T\wedge\bar\partial_t\Gamma_t^{n-2,\,n-1} & = & -\int\limits_{X_t}T\wedge(\bar\partial_t\overline{\Gamma_t^{n-2,\,n-1}} + \partial_t\Gamma_t^{n-2,\,n-1})\\
   & = & -\int\limits_{X_t}T\wedge(\gamma_t^{n-1} - \Omega_{A,\,t}^{n-1,\,n-1}).\end{eqnarray}

\noindent Now, $\int_{X_t}T\wedge\gamma_t^{n-1}$ converges to $\int_{X_0}T\wedge\gamma_0^{n-1}$ and, by (\ref{eqn:Omega-n-1-conv}), $\int_{X_t}T\wedge\Omega_{A,\,t}^{n-1,\,n-1}$ converges to $\int_{X_0}T\wedge\Omega_{A,\,0}'$ as $t\rightarrow 0$. Putting together (\ref{eqn:1main-integral_t}), (\ref{eqn:2main-integral_t}) and all the pieces of convergence information just mentioned, we get the convergence\!\!:

\begin{equation}\label{eqn:3main-integral_t}\int\limits_{X_t}T\wedge(\Omega_t^{n,\,n-2} + \gamma_t^{n-1} + \overline{\Omega_t^{n,\,n-2}})\longrightarrow \int\limits_{X_0}T\wedge\Omega_{A,\,0}' = \int\limits_{X_0}T\wedge\gamma_0^{n-1}    \hspace{2ex} \mbox{as}\hspace{1ex}t\rightarrow 0,\end{equation}

\noindent where the last identity is nothing but (\ref{eqn:T_to-show-bis}).

 Recall that $\int_{X_t}T\wedge(\Omega_t^{n,\,n-2} + \gamma_t^{n-1} + \overline{\Omega_t^{n,\,n-2}})\geq 0$ for all $t\sim 0$ with $t\neq 0$ by (\ref{eqn:integral-choiceT}). Hence (\ref{eqn:3main-integral_t}) implies $\int\limits_{X_0}T\wedge\gamma_0^{n-1}\geq 0$ and we are done. \hfill $\Box$

\section{Relations between the sGG class and other classes of compact complex manifolds}\label{section:sGG-other}

In this section we show that sGG manifolds are unrelated to balanced manifolds
and to those whose Fr\"olicher spectral sequence degenerates at $E_1$. Examples of compact complex manifolds $X$ with ${\cal SG}_X \not= {\cal G}_X$ but admitting strongly Gauduchon metrics are also given. To construct appropriate examples we will consider the class of nilmanifolds endowed with an invariant complex structure.

Recall that a nilmanifold $N=G/\Gamma$ is a compact quotient of a connected and simply-connected nilpotent real Lie group $G$
by a lattice $\Gamma$ of maximal rank in $G$.
Let $\mathfrak{g}$ be the Lie algebra of the group $G$. We will say that \emph{``$N$ has underlying Lie algebra $\mathfrak{g}$''}
or that \emph{``$\mathfrak{g}$ is the Lie algebra underlying $N$''}.
We will denote 6-dimensional real Lie algebras in the usual abbreviated form; for instance,
$(0^4,12,34)$ denotes the Lie algebra $\mathfrak{g}$ with generators $\{e_i\}_{i=1}^6$ satisfying
the bracket relations $[e_1,e_2]=-e_5,\ [e_3,e_4]=-e_6$, or equivalently there exists a basis $\{\alpha^i\}_{i=1}^6$
of the dual $\mathfrak{g}^*$ such that
$d\alpha^1=d\alpha^2=d\alpha^3=d\alpha^4=0,\ d\alpha^5=\alpha^1\wedge \alpha^2,\ d\alpha^6=\alpha^3\wedge \alpha^4$.

Notice that by Nomizu's theorem [Nom54], the integer $k$ appearing in $0^k$ in the notation above
is precisely the first Betti number of $N$, i.e. $b_1(N)=k$.

The complex structures that we will consider on $N$ are invariant in the sense that
they stem naturally from ``complex'' structures $J$ on the Lie algebra $\mathfrak{g}$ of $G$.
For any such $J$, the $i$-eigenspace $\mathfrak{g}_{1,0}$ of $J$ in
$\mathfrak{g}_\mathbb{C}=\mathfrak{g} \otimes_\mathbb{R} \mathbb{C}$
is a complex subalgebra. When $\mathfrak{g}_{1,0}$ is abelian we will refer to $J$ as
an {\it abelian} complex structure.

The following result identifies the compact complex nilmanifolds of complex dimension 3
that are sGG.


\begin{The}\label{sGGnilm}
Let $N$ be a nilmanifold of (real) dimension six not isomorphic to a torus and let $J$ be an invariant complex structure on $N$. Then, the compact complex manifold $X=(N,J)$ is sGG if and only if the Lie algebra underlying $N$ is isomorphic to
$(0^4,12,34)$, $(0^4,12,14+23)$, $(0^4,13+42,14+23)$ or $(0^4,12,13)$ and the complex structure $J$ is not abelian.
\end{The}

\noindent {\it Proof.}
By Theorem~\ref{The:Betti-dbar}, if $N$ admits an invariant complex structure $J$ such that $X=(N,J)$ is sGG then the first Betti number is even.
From the classification of nilpotent Lie algebras admitting a complex structure [Sal01], this condition implies that the Lie algebra underlying $N$ belongs to the following list\!\!:
$(0^4,12,34)$, $(0^4,12,14+23)$, $(0^4,13+42,14+23)$, $(0^4,12,13)$,
$(0^4,12,14+25)$, $(0^2,12,13,23,14+25)$.

We first rule out the last two cases. It was proved in [UV14, Proposition 2.4] that for any invariant complex structure $J$ on a nilmanifold $N$ with underlying Lie algebra $(0^2,12,13,23,14+25)$
there is a global basis $\{\eta^1,\eta^2,\eta^3\}$ of forms of bidegree (1,0), with respect to $J$, satisfying complex equations of the shape\!\!:
$$
d\eta^1=0,\ \ d\eta^2=\eta^{13} + \eta^{1\bar3},\ \
d\eta^3=i\eta^{1\bar 1} \pm i(\eta^{1\bar2} -\eta^{2\bar1}),
$$
(where we use the standard notation\!\!: $\eta^{jk}:=\eta^j\wedge\eta^k$, $\eta^{j\bar{k}}:=\eta^j\wedge\overline{\eta^k}$.)

\noindent That is to say, up to equivalence there exist exactly two invariant complex structures on $N$ depending on the choice of sign in the third equation. Hence,
$H^{0,\,1}_{\bar\partial}(N,J)=\langle [\eta^{\bar 1}]_{\bar\partial}, [\eta^{\bar 3}]_{\bar\partial} \rangle$ by a result in [Rol09, Section 4.2],
and we get $b_1(N)=2<4= 2 h^{0,\,1}_{\bar\partial}(N,J)$.
It follows from Theorem~\ref{The:Betti-dbar} that there is no invariant complex structure on $N$ satisfying the sGG property.

By [COUV11],
for any invariant complex structure $J$ on a nilmanifold $N$ with underlying Lie algebra $(0^4,12,14+25)$
there is a global basis $\{\eta^1,\eta^2,\eta^3\}$ of (1,0)-forms satisfying
$$
d\eta^1=0,\ \ d\eta^2=\eta^{1\bar1},\ \ d\eta^3=\eta^{1\bar2} +\eta^{2\bar1},
$$
which implies that $h^{0,\,1}_{\bar\partial}(N,J)=3$. Therefore, $b_1(N)=4<6= 2 h^{0,\,1}_{\bar\partial}(N,J)$, so there is no invariant complex structure on $N$ satisfying the sGG property.

It is well known that on $6$-dimensional nilmanifolds different from the complex tori there exists (up to equivalence) only one complex-parallelisable complex
structure given by the complex equations
$$
d\eta^1=d\eta^2=0,\quad d\eta^3=\eta^{12}.
$$
This corresponds to the Iwasawa manifold (which is an sGG manifold by e.g. Corollary \ref{Cor:Iwasawa-sGG}) whose underlying Lie algebra is precisely $(0^4,13+42,14+23)$.
Now, for any other complex structure $J$ on a nilmanifold $N$ with underlying Lie algebra $(0^4,12,34)$, $(0^4,12,14+23)$, $(0^4,13+42,14+23)$ or $(0^4,12,13)$, it is proved in [COUV11] that there is a (1,0)-basis satisfying
\begin{equation}\label{eqn:complex-equations}
d\eta^1=d\eta^2=0,\quad d\eta^3=\rho\,\eta^{12} + \eta^{1\bar1} + \lambda\,\eta^{1\bar2} + D\,\eta^{2\bar2},
\end{equation}
where $\rho\in\{0,1\}$, $\lambda\in\mathbb R^{\geq 0}$ and $D=x+iy\in\mathbb C$ with $y\geq 0$.
Notice that $J$ is abelian if and only if $\rho=0$.

To complete the proof we must show that for any complex structure $J$ given by \eqref{eqn:complex-equations}
the compact complex manifold $(N,J)$ satisfies $b_1(N)=4= 2 h^{0,\,1}_{\bar\partial}(N,J)$ if and only if $\rho=1$.
But this is clear because
$H^{0,\,1}_{\bar\partial}(N,J)=\langle [\eta^{\bar 1}]_{\bar\partial}, [\eta^{\bar 2}]_{\bar\partial}, [\eta^{\bar 3}]_{\bar\partial} \rangle$
when $\rho=0$,
and
$H^{0,\,1}_{\bar\partial}(N,J)=\langle [\eta^{\bar 1}]_{\bar\partial}, [\eta^{\bar 2}]_{\bar\partial} \rangle$
for $\rho=1$.
\hfill $\Box$

\vspace{2ex}

For any compact complex manifold $X$, it is immediate that if $X$ is sGG then
$X$ has an sG metric. The following example shows that
the converse does not hold in general even if the sG hypothesis is reinforced to the balanced hypothesis and even with an extra property.

\begin{Prop}\label{no-converse}
There exists a compact complex manifold $X$ having a balanced metric, with Fr\"olicher spectral sequence
degenerating at the first step and with first Betti number $b_1(X)$ odd. Thus, $X$ is not sGG.
\end{Prop}

\noindent {\it Proof.}
Let $N$ be a nilmanifold with underlying Lie algebra isomorphic to $(0^5,12+34)$. Then, the first Betti number of $N$
is 5.
We consider on $N$ the complex structure $J$ defined by
the complex equations
$$
d\eta^1=d\eta^2=0,\ \ d\eta^3=\eta^{1\bar{1}}-\eta^{2\bar{2}}.
$$
By Theorem~\ref{The:Betti-dbar} we know that $X=(N,J)$ is not sGG because $b_1(N)=5$.
It is proved in [COUV11] that $E_1(X)\cong E_\infty(X)$.
Moreover, $X$ is balanced;
for instance,
$\omega=\frac{i}{2}(\eta^{1\bar{1}}+\eta^{2\bar{2}}+\eta^{3\bar{3}})$
satisfies $d\omega^2=0$, that is, $\omega$ is a balanced metric on $X$.
\hfill $\Box$

\begin{Prop}\label{unrelated-sGG-balanced-Frolicher}
The balanced property and the sGG property are unrelated.
Moreover, the Fr\"olicher spectral sequence degenerating at $E_1$ and the sGG property are also unrelated.
\end{Prop}

\noindent {\it Proof.}
In Proposition~\ref{no-converse} we proved that ``balanced'' does not imply ``sGG'',
and that $E_1(X)\cong E_\infty(X)$ does not imply $X$ to be sGG.
We now show that there exists an sGG compact complex manifold $X$ that is not balanced and whose Fr\"olicher spectral sequence does not degenerate at $E_1$.

Let $N$ be a nilmanifold with underlying Lie algebra isomorphic to $(0^4,13+42,14+23)$, that is, $N$ is the (real)
manifold underlying the Iwasawa manifold. We consider on $N$ the complex structure $J$ defined by
the complex equations
\begin{equation}\label{example1}
d\eta^1=d\eta^2=0,\ \ d\eta^3=\eta^{12}+\eta^{1\bar{1}}.
\end{equation}
By Theorem~\ref{sGGnilm} the compact complex manifold $X=(N,J)$ is sGG because the complex structure $J$ is not abelian.
However, from the general study in [COUV11] one has that $E_1(X)\not\cong E_2(X) \cong E_\infty(X)$
and $X$ does not admit any balanced metric.
\hfill $\Box$

\vspace{3ex}

We are now ready to give another example promised in the introduction.

\begin{Prop}\label{unrelated-sGG-super-sG}
The superstrong Gauduchon property and the sGG property are unrelated.
\end{Prop}

\noindent {\it Proof.}
Let us show first that there exists an sGG compact complex manifold $X$ that does not admit any superstrong Gauduchon metric.
Consider $X=(N,J)$ a (real) $2n$-dimensional nilmanifold $N$ endowed with an invariant complex structure $J$.
By the usual symmetrisation process, if $\omega$ is a superstrong Gauduchon metric on $X$, then there also exists an invariant
superstrong Gauduchon metric $\hat\omega$ on $X$. Indeed, if $\Omega=\omega^{n-1}$ satisfies $\partial\Omega=\partial\bar\partial \alpha$ for some
$(n-1,n-2)$-form $\alpha$, then by symmetrisation we get that the positive definite invariant $(n-1,n-1)$-form $\widetilde\Omega$ (obtained from $\Omega$) satisfies
$\partial\widetilde\Omega=\partial\bar\partial\tilde\alpha$ for an invariant $(n-1,n-2)$-form $\tilde\alpha$. Now, since $\widetilde\Omega>0$, it is well known that
there exists an invariant Hermitian metric $\hat\omega$ such that $\widetilde\Omega=\hat\omega^{n-1}$. Thus $\hat\omega$ is necessarily an invariant superstrong Gauduchon metric on $X$.

Now, let us consider $X=(N,J)$ defined by \eqref{example1}, which by Theorem~\ref{sGGnilm} is sGG.
A direct calculation shows that
$\partial\bar\partial \Lambda^{2,1} (\mathfrak{g}^*) \equiv 0$. Therefore, if a superstrong Gauduchon metric existed on $X$, it would have to be an invariant balanced metric. However, we pointed out in the proof of
Proposition~\ref{unrelated-sGG-balanced-Frolicher} that $X$ is not balanced. Thus, $X$ is sGG but does not admit any superstrong Gauduchon metric.

Conversely, we notice that the superstrong Gauduchon property does not imply the sGG property because, thanks to Proposition~\ref{no-converse}, there exists a balanced manifold which is not sGG.   \hfill $\Box$

\section{Examples of deformation limits of sGG manifolds}\label{sGG-nonclosedness}

The following result shows that the sGG hypothesis on $X$ does not ensure the Bott-Chern number $h^{1,\,1}_{BC}(X)$ to be locally deformation constant.

\begin{Prop}\label{Prop:h11-jumping}
There exists a holomorphic family of compact complex {\bf sGG} manifolds $(X_t)_{t\in\Delta}$ such that $h^{1,\,1}_{BC}(0) > h^{1,\,1}_{BC}(t)$
for all $t\in\Delta\setminus C$, where $\Delta\subset\C$ is a small open disc about $0$ and $C$ is a real curve through $0$.
\end{Prop}

\noindent {\it Proof.}
Let $X_0=(N,J_0)$ be a complex nilmanifold of real dimension $6$ defined by the equations
\begin{equation}\label{example22}
d\eta^1=d\eta^2=0,\quad d\eta^3=\eta^{12}+\eta^{1\bar{1}}+\eta^{1\bar{2}}-2\,\eta^{2\bar{2}}.
\end{equation}
By [COUV11, Table 1] the Lie algebra $\mathfrak{g}$ underlying $N$ is isomorphic to $(0^4,12,14+23)$.
Since the complex structure $J_0$ is not abelian, Theorem~\ref{sGGnilm} implies that $X_0$ is sGG.

By [Ang11, Theorem 2.7], the Bott-Chern cohomology groups of $X_0$
can be calculated at the level of the Lie algebra underlying $N$, in particular,
$H^{1,\,1}_{BC}(X_0)\cong H^{1,\,1}_{BC}(\mathfrak{g},J_0)=\ker\{d\colon \Lambda^{1,1}(\mathfrak{g}^*)\longrightarrow \Lambda^3(\mathfrak{g}^*_\mathbb{C})\}$.
From the equations \eqref{example22} we get
$$
H^{1,\,1}_{BC}(X_0) \cong \langle [\eta^{1\bar{1}}]_{BC}, [\eta^{1\bar{2}}]_{BC}, [\eta^{2\bar{1}}]_{BC}, [\eta^{2\bar{2}}]_{BC},
[\eta^{1\bar{3}}\!+\!2\,\eta^{2\bar{3}}\!+\!\eta^{3\bar{1}}\!+\!2\,\eta^{3\bar{2}}]_{BC} \rangle,
$$
therefore $h^{1,\,1}_{BC}(X_0)=5$.

Now we consider a small deformation $J_t$ given by
$$
t \frac{\partial}{\partial z_2} \otimes d \bar{z}_2 \in H^{0,\,1}(X_0,\,T^{1,0}X_0),
$$
where $z_2$ is a complex coordinate such that $\eta^2=d z_2$.
By Corollary~\ref{Cor:openness} we know that the compact complex manifold $X_t=(N,J_t)$ is sGG for all $t \in \C$ close enough to $0$.
In fact, for $t\in \mathbb{C}$ with $|t|<1$, if we consider the basis
$\{\nu^1_t=\eta^1, \nu^2_t=\frac{1-\bar{t}}{1-|t|^2}(\eta^2+t\,\eta^{\bar{2}}), \nu^3_t=\eta^3\}$
of complex forms of type (1,0) with respect to $J_t$, then
the complex structure equations along the deformation are:
\begin{equation}\label{example22bis}
d\nu^1_t=d\nu^2_t=0,\quad
d\nu^3_t=\nu_t^{12}+\nu_t^{1\bar{1}}+\nu_t^{1\bar{2}}-2\,\frac{1-|t|^2}{|1-t|^2}\,\nu_t^{2\bar{2}}.
\end{equation}

Next we compute the dimension of the Bott-Chern cohomology group $H^{1,\,1}_{BC}(X_t)$ of $X_t$.
Since the complex structure $J_t$ is invariant, we can use again [Ang11, Theorem 2.7] to reduce the calculation
to the invariant forms.
By \eqref{example22bis} it is clear that $\nu_t^{1\bar{1}}$, $\nu_t^{1\bar{2}}$, $\nu_t^{2\bar{1}}$ and $\nu_t^{2\bar{2}}$
define Bott-Chern classes in $H^{1,\,1}_{BC}(X_t)$. To see if there are some other classes, we need to compute the differentials
of the remaining basic (1,1)-forms $\nu_t^{j\bar{k}}$. From the equations \eqref{example22bis} we get:

\vspace{2ex}

$d\nu_t^{1\bar{3}}=\nu_t^{12\bar{1}}-2\,\frac{1-|t|^2}{|1-t|^2}\nu_t^{12\bar{2}}-\nu_t^{1\bar{1}\bar{2}}$,

\vspace{2ex}

$d\nu_t^{2\bar{3}}=-\nu_t^{12\bar{1}}-\nu_t^{2\bar{1}\bar{2}},$

\vspace{2ex}

$d\nu_t^{3\bar{1}}=\nu_t^{12\bar{1}}-\nu_t^{1\bar{1}\bar{2}}+2\,\frac{1-|t|^2}{|1-t|^2}\nu_t^{2\bar{1}\bar{2}}$,

\vspace{2ex}

$d\nu_t^{3\bar{2}}=\nu_t^{12\bar{2}}+\nu_t^{1\bar{1}\bar{2}},$

\vspace{2ex}

$d\nu_t^{3\bar{3}}=\nu_t^{12\bar{3}}-\nu_t^{13\bar{1}}-\nu_t^{23\bar{1}}+2\,\frac{1-|t|^2}{|1-t|^2}\nu_t^{23\bar{2}} +\nu_t^{1\bar{1}\bar{3}}+\nu_t^{1\bar{2}\bar{3}}
-2\,\frac{1-|t|^2}{|1-t|^2}\nu_t^{2\bar{2}\bar{3}}-\nu_t^{3\bar{1}\bar{2}}.$

\vspace{2ex}

\noindent From these expressions, it is easy to check that there exists at most one more closed (1,1)-form, and that such a form exists if and only if
$1-|t|^2=|1-t|^2$.

Let $C=\{t\in\C \mid |t|^2+|1-t|^2=1 \}$. Note that $C$ is a circle centered at $t=1/2$ passing through $t=0$.
Our discussion above shows that $h^{1,\,1}_{BC}(X_t)=4$
for all $t\in\Delta^{\star}\setminus C$, where $\Delta=\{t\in\C \mid |t|<1\}\subset\C$,
that is, the Bott-Chern number $h^{1,\,1}_{BC}$ is not locally deformation constant.
\hfill $\Box$

\vspace{2ex}

In the following result we show by means of three examples that the sGG property of compact complex manifolds is not closed under holomorphic deformations. The behaviour of the holomorphic families in the three examples is different and illustrate several
possibilities for the central limit.

\begin{Prop}\label{not-closed}
There exist holomorphic families of compact complex manifolds $(X_t)_{t\in\Delta}$ over an open disc $\Delta\subset\C$ about $0$ such that $X_t$ is sGG for all $t\in\Delta\setminus\{0\}$, but $X_0$ is not sGG.
\end{Prop}

\noindent {\it Proof.} We will describe three examples in succession.

\noindent {\it First example.} Let us consider the compact complex manifold $X_0=(N,J_0)$, where
$N$ is the nilmanifold with underlying Lie algebra $(0^4,13+42,14+23)$ and $J_0$ is the abelian structure defined by
the complex structure equations
$$
d\eta^1=d\eta^2=0,\quad d\eta^3=\eta^{1\bar{1}}+\eta^{1\bar{2}}.
$$
Since $J_0$ is abelian, the manifold $X_0$ is not sGG (by Theorem~\ref{sGGnilm}), and
$H^{0,1}_{\bar\partial}(X_0,\C)=\langle [\eta^{\bar{1}}]_{\bar\partial}, [\eta^{\bar{2}}]_{\bar\partial}, [\eta^{\bar{3}}]_{\bar\partial} \rangle$.
Consider a small deformation $J_t$ given by
$$
t \frac{\partial}{\partial z_2} \otimes d \bar{z}_2 \in H^{0,\,1}(X_0,\,T^{1,0}X_0),
$$
where $z_2$ is a complex coordinate such that $\eta^2=d z_2$.
Let us consider the basis
$\{\tau^1_t=\eta^1, \tau^2_t=\eta^2+t\,\eta^{\bar{2}}, \tau^3_t=\eta^3\}$
of complex forms of type (1,0) with respect to $J_t$.
Then, for $t\in \mathbb{C}$ with $|t|<1$, the complex structure equations of the deformation are:
$$
d\tau^1_t=d\tau^2_t=0,\quad
d\tau^3_t=-\frac{\bar{t}}{1-|t|^2}\, \tau_t^{12}+\tau_t^{1\bar{1}}+\frac{1}{1-|t|^2}\, \tau_t^{1\bar{2}}.
$$
For any $t\not=0$, the complex structure is not abelian because the differential of the (1,0)-form $\tau^3_t$
has a non-zero component of bidegree (2,0), so the compact complex manifold $X_t=(N,J_t)$ is sGG for any
$t\not=0$ by Theorem~\ref{sGGnilm}.

\vspace{1ex}

Note that $X_0$ admits a balanced metric by [COUV11, Proposition 7.7], hence also an sG metric, so $\mathcal{SG}_{X_0}\not=\emptyset$.

\vspace{2ex}

\noindent {\it Second example.} In [COUV11, Theorem 7.9] it is constructed a holomorphic family
of compact complex (nil)manifolds $(X_t)_{t\in\Delta}$ over an open disc $\Delta\subset\C$ about $0$, where $X_t$ is balanced for any $t\not=0$
and such that the central
limit $X_0$ is a complex nilmanifold with underlying Lie algebra $(0^4,12,14+23)$ endowed with an abelian complex structure $J_0$.
The complex structure on $X_t$ is invariant and non-abelian for any $t\not=0$, so by Theorem~\ref{sGGnilm} the compact complex manifold
$X_t$ is sGG, but the central limit $X_0$ is not sGG. Moreover, it is proved in [COUV11, Proposition 7.7] that $X_0$ is not sG,
so $\mathcal{SG}_{X_0}=\emptyset$.

\vspace{2ex}

\noindent {\it Third example.} Angella and Kasuya obtain in [AK14, Proposition 4.1 (i)] a holomorphic family
of compact complex manifolds $X_t$ over an open disc in $\C$ about $0$, satisfying the $\partial\bar\partial$-lemma for any $t\not=0$
and such that the central
limit $X_0$ is the complex-parallelisable Nakamura manifold [Nak75]. By Theorem~\ref{The:Betti-dbar} we conclude that $X_0$
is not sGG because $b_1(X_0)=2<6=2 h^{0,\,1}_{\bar\partial}(X_0)$ (see [AK14, Table 10]).
Note however that the central limit $X_0$ is balanced.
\hfill $\Box$

\vspace{2ex}

In [Pop09, Proposition 4.1] it is proved that given a holomorphic family of compact complex manifolds $(X_t)_{t\in\Delta}$ over an open disc $\Delta\subset\C$ about $0$, if $X_t$ satisfies the $\partial\bar\partial$-lemma for all $t\in\Delta\setminus\{0\}$ then $X_0$ is sG.
However, the central limit $X_0$ may be neither sGG (see {\it Third example} in the proof of Proposition~\ref{not-closed})
nor balanced (see [FOU14, Theorem 5.2]).
Furthermore, in the following proposition we show that in general $X_0$ does not admit superstrong Gauduchon metrics.

\begin{Prop}\label{central-limit-not-super-sG}
There exists a holomorphic family of compact complex manifolds $(X_t)_{t\in\Delta}$ over an open disc $\Delta\subset\C$ about $0$ such that $X_t$
satisfies the $\partial\bar\partial$-lemma for all $t\in\Delta\setminus\{0\}$, but $X_0$ is not superstrong Gauduchon.
\end{Prop}

\noindent {\it Proof.}
We consider the holomorphic family of compact complex manifolds $(X_t)_{t\in\Delta}$ constructed in [FOU14, Theorem 5.2], which satisfies the
$\partial\bar\partial$-lemma for all $t\in\Delta\setminus\{0\}$. The central limit of that family is $X_0=(G/\Gamma,J_0)$, where $G/\Gamma$ is a solvmanifold (i.e. a compact quotient of a connected and simply-connected {\it solvable} real Lie group $G$
by a lattice $\Gamma$ of maximal rank in $G$) and $J_0$ is the invariant complex structure defined by the complex structure equations
$$
d\eta^1=2i\,\eta^{13}+\eta^{3\bar{3}},\quad
d\eta^2=-2i\,\eta^{23},\quad
d\eta^3=0,
$$
where $\{\eta^1,\eta^2,\eta^3\}$ is a (1,0)-basis. Thus, we can apply the symmetrisation process
and proceed as in the proof of Proposition~\ref{unrelated-sGG-super-sG}.
So it suffices to show that there do not exist superstrong Gauduchon metrics on the underlying solvable Lie algebra $\mathfrak{g}$.
From the complex structure equations above, it is easy to check that
$\partial\bar\partial \Lambda^{2,1} (\mathfrak{g}^*) \equiv 0$, which implies that
any invariant superstrong Gauduchon metric must be balanced. But this is not possible by [FOU14, Theorem 5.2],
so we conclude that $X_0$ is not superstrong Gauduchon. (However, $X_0$ is sG as pointed out in [FOU14].)
\hfill $\Box$

\vspace{3ex}

\noindent {\bf Acknowledgments.}
This work has been partially supported through Project MICINN (Spain) MTM2011-28326-C02-01.

\vspace{3ex}

\noindent {\bf References} \\

\noindent [AB95]\, L. Alessandrini, G. Bassanelli --- {\it Modifications of Compact Balanced Manifolds} --- C. R. Acad. Sci. Paris, S\'er. I {\bf 320} (1995), 1517-1522.

\vspace{1ex}

\noindent [Ang11]\, D. Angella --- {\it The Cohomologies of the Iwasawa Manifold and of Its Small Deformations} --- J. Geom. Anal. (2011) DOI: 10.1007/s12220-011-9291-z.

\vspace{1ex}

\noindent [ADT14]\, D. Angella, G. Dloussky, A. Tomassini --- {\it On Bott-Chern Cohomology of Compact Complex Surfaces} --- Annali di Matematica Pura ed Applicata, doi 10.1007/s10231-014-0458-7

\vspace{1ex}

\noindent [AK12]\, D. Angella, H. Kasuya, --- {\it Bott-Chern cohomology of solvmanifolds} --- arXiv e-print DG 1212.5708v3.

\vspace{1ex}

\noindent [AK14]\, D. Angella, H. Kasuya, --- {\it Cohomologies of deformations of solvmanifolds and closedness of some properties} --- arXiv e-print CV 1305.6709v2, to appear in Mathematica Universalis.

\vspace{1ex}

\noindent [BDPP]\, S. Boucksom, J.-P. Demailly, M. Paun, T. Peternell --- {\it The Pseudo-effective Cone of a Compact K\"ahler Manifold and Varieties of Negative Kodaira Dimension} --- J. Alg. Geom. {\bf 22} (2013) 201-248.

\vspace{1ex}

\noindent [Buc99]\, N. Buchdahl --- {\it On Compact K\"ahler Surfaces} --- Ann. Inst. Fourier {\bf 49}, no. 1 (1999) 287-302.

\vspace{1ex}

\noindent [COUV11]\,  M. Ceballos, A. Otal, L. Ugarte, R. Villacampa --- {\it Invariant Complex Structures on 6-nilmanifolds: Classification, Fr\"olicher Spectral Sequence and Special  Hermitian Metrics} --- J. Geom. Anal., doi 10.1007/s12220-014-9548-4

\vspace{1ex}

\noindent [Dem92]\, J.-P. Demailly --- {\it Regularization of Closed Positive Currents and Intersection Theory} --- J. Alg. Geom., {\bf 1} (1992), 361-409.

\vspace{1ex}

\noindent [DP04]\, J.-P. Demailly, M. Paun --- {\it Numerical Charaterization of the K\"ahler Cone of a Compact K\"ahler Manifold} --- Ann. Math. (2) {\bf 159(3)} (2004) 1247-1274.

\vspace{1ex}

\noindent [FOU14]\, A. Fino, A. Otal, L. Ugarte --- {\it Six Dimensional Solvmanifolds with Holomorphically Trivial Canonical Bundle} --- Int. Math. Res. Not. IMRN, doi 10.1093/imrn/rnv112

\vspace{1ex}

\noindent [Fuj78]\, A. Fujiki --- {\it Closedness of the Douady Spaces of Compact K\"ahler Spaces} --- Publ. RIMS, Kyoto Univ. {\bf 14} (1978), 1-52.

\vspace{1ex}

\noindent [Gau77a]\, P. Gauduchon --- {\it Le th\'eor\`eme de l'excentricit\'e nulle} --- C.R. Acad. Sc. Paris, S\'erie A, t. {\bf 285} (1977), 387-390.

\vspace{1ex}

\noindent [Gau77b]\, P. Gauduchon --- {\it Fibr\'es hermitiens \`a endomorphisme de Ricci non n\'egatif} --- Bull. Soc. Math. France {\bf 105} (1977) 113-140.

\vspace{1ex}

\noindent [Har77]\, R. Hartshorne --- {\it Algebraic Geometry} --- Springer, Graduate Texts in Mathematics, {\bf 52} (1977).

\vspace{1ex}

\noindent [KS60]\, K. Kodaira, D.C. Spencer --- {\it On Deformations of Complex Analytic Structures, III. Stability Theorems for Complex Structures} --- Ann. Math. {\bf 71}, No. 1 (1960), 43-76.

\vspace{1ex}

\noindent [Lam99]\, A. Lamari --- {\it Courants k\"ahl\'eriens et surfaces compactes} --- Ann. Inst. Fourier, Grenoble, {\bf 49}, 1 (1999), 263-285.

\vspace{1ex}

\noindent [Mic83]\, M. L. Michelsohn --- {\it On the Existence of Special Metrics in Complex Geometry} --- Acta Math. {\bf 143} (1983) 261-295.



\vspace{1ex}

\noindent [Miy74]\, Y. Miyaoka --- {\it K\"ahler Metrics on Elliptic Surfaces} --- Proc. Japan Acad. {\bf 50} No. 8 (1974) 533-536.

\vspace{1ex}

\noindent [Nak75]\, I. Nakamura --- {\it Complex Parallelisable Manifolds and Their Small Deformations} --- J. Diff. Geom. {\bf 10} (1975), 85-112.

\vspace{1ex}

\noindent [Nom54] \, K. Nomizu --- {\it On the cohomology of compact homogeneous spaces of nilpotent Lie groups}--- Ann. Math. {\bf 59} (1954), 531-538.

\vspace{1ex}

\noindent [Pop09]\, D. Popovici --- {\it Limits of projective manifolds under holomorphic deformations} --- arXiv e-print AG 0910.2032v1

\vspace{1ex}

\noindent [Pop13a]\, D. Popovici --- {\it Deformation Limits of Projective Manifolds\!\!: Hodge Numbers and Strongly Gauduchon Metrics} --- Invent. Math. {\bf 194} (2013), 515-534.

\vspace{1ex}

\noindent [Pop13b]\, D. Popovici --- {\it Aeppli Cohomology Classes Associated with Gauduchon Metrics on Compact Complex Manifolds} --- Bull. Soc. Math. France {\bf 143} (3), (2015), p. 1-37.

\vspace{1ex}

\noindent [Pop14]\, D. Popovici --- {\it Deformation Openness and Closedness of Various Classes of Compact Complex Manifolds; Examples} --- Ann. Sc. Norm. Super. Pisa Cl. Sci. (5), Vol. XIII (2014), 255-305.

\vspace{1ex}

\noindent [Rol09]\, S. Rollenske --- {\it Geometry of nilmanifolds with left-invariant
complex structure and deformations in the large} --- Proc. London Math. Soc. {\bf 99} (2009), 425-460.

\vspace{1ex}

\noindent [Sal01]\, S. Salamon --- {\it Complex structures on nilpotent Lie algebras} ---
J. Pure Appl. Algebra {\bf 157} (2001), 311-333.

\vspace{1ex}

\noindent [Sch07]\, M. Schweitzer --- {\it Autour de la cohomologie de Bott-Chern} --- arXiv e-print math.AG/0709.3528v1.

\vspace{1ex}

\noindent [Siu83]\, Y.-T. Siu --- {\it Every K3 Surface Is K\"ahler} --- Invent. Math. {\bf 73} (1983) 139-150.



\vspace{1ex}

\noindent [UV14]\, L. Ugarte, R. Villacampa --- {\it Non-nilpotent complex geometry of nilmanifolds and heterotic
supersymmetry} --- Asian J. Math. {\bf 18} (2014), 229-246.

\vspace{1ex}

 \noindent [Var86]\, J. Varouchas --- {\it Sur l'image d'une vari\'et\'e k\"ahl\'erienne compacte} --- LNM {\bf 1188}, Springer (1986) 245--259.

\vspace{3ex}

\noindent Institut de Math\'ematiques de Toulouse, Universit\'e Paul Sabatier,

\noindent 118 route de Narbonne, 31062 Toulouse, France

\noindent Email\!: popovici@math.univ-toulouse.fr

\vspace{2ex}

\noindent and

\vspace{2ex}

\noindent Departamento de Matem\'aticas\,-\,I.U.M.A., Universidad de Zaragoza,

\noindent Campus Plaza San Francisco, 50009 Zaragoza, Spain

\noindent Email\!: ugarte@unizar.es

\end{document}